\documentclass[a4paper,DIV=12,11pt,parskip=half]{scrartcl}

\usepackage[USenglish]{babel}
\usepackage[utf8x]{inputenc}
\usepackage[pdfencoding=auto,
unicode=true,
plainpages=false,
pdfpagelabels=true]{hyperref}

\usepackage{textcomp}
\newcommand{\qdist}[1]{\ifmmode\langle#1\rangle\else\textlangle#1\textrangle\fi}

\usepackage{amssymb,amsmath}

\usepackage{mathtools}
\mathtoolsset{centercolon}

\usepackage{graphicx}
\graphicspath{{graphics/}}

\usepackage{multirow}
\usepackage{booktabs}

\usepackage[binary-units=true]{siunitx}
\usepackage[finalizecache=false,
frozencache=true]{minted}

\usepackage{caption}
\usepackage{subcaption}
\usepackage{amsthm}

\usepackage{enumitem}

\usepackage[ruled,vlined]{algorithm2e}

\usepackage[table]{xcolor}
\usepackage[section]{placeins}

\usepackage{cleveref}
\crefname{section}{Sec.}{Sec.}
\Crefname{equation}{Eq.}{Eqs.}
\Crefname{figure}{Fig.}{Figs.}
\Crefname{tabular}{Tab.}{Tabs.}

\renewcommand{\vec}{\boldsymbol}
\newcommand{\mat}{\boldsymbol}

\renewcommand{\d}{\; \mathrm d}
\renewcommand{\P}{\mathbb{P}}
\newcommand{\R}{\mathbb{R}}
\newcommand{\N}{\mathbb{N}}

\newcommand{\textblue}[1]{{#1}}
\newcolumntype{b}{>{}c}

\theoremstyle{definition}

\newtheorem{problem}{Problem}
\newtheorem{remark}{Remark}

\theoremstyle{remark}

\numberwithin{equation}{section}
\numberwithin{figure}{section}
\numberwithin{table}{section}
\numberwithin{problem}{section}
\numberwithin{remark}{section}

\begin{document}



\title{A geometric multigrid method for space-time finite element discretizations of the Navier--Stokes equations and its application to 3d flow simulation}


\author{Mathias Anselmann$^{\star,1}$, Markus Bause$^{\star}$\\
	{\small ${}^\star$ Helmut Schmidt University, Faculty of 
		Mechanical Engineering, Holstenhofweg 85,}\\ 
	{\small 22043 Hamburg, Germany}\\
}

\date{}
\maketitle

\footnotetext[1]{Corresponding author: anselmann@hsu-hh.de}

\begin{abstract}
\textbf{Abstract.}
We present a parallelized geometric multigrid (GMG) method, based on the cell-based Vanka smoother, for higher order space-time finite element methods (STFEM) to the incompressible Navier--Stokes equations. The STFEM is implemented as a time marching scheme. The GMG solver is applied as a preconditioner for GMRES iterations. Its performance properties are demonstrated for 2d and 3d benchmarks of flow around a cylinder. The key ingredients of the GMG approach are the construction of the local Vanka smoother over all degrees of freedom in time of the respective subinterval and its efficient application. For this, data structures that store pre-computed cell inverses of the Jacobian for all hierarchical levels and require only a reasonable amount of memory overhead are generated. The GMG method is built for the \emph{deal.II} finite element library. The concepts are flexible and can be transferred to similar  software platforms.
\end{abstract}

\section{Introduction}

The accurate numerical simulation of incompressible viscous flow continues to remain a challenging task, in particular, if three space dimensions are considered due to the sake of physical realism. Higher order methods offer the potential to achieve accurate results on computationally feasible grids with a minimum of numerical costs. However, constructing higher order numerical methods maintaining stability and inheriting most of the rich structure of the continuous problem becomes an important prerequisite. For this, we refer, e.g., to \cite{johnHigherorderFiniteElement2001,johnHigherOrderFinite2002} for stationary problems and to \cite{hussainHigherOrderGalerkin2013,hussainEfficientStableFinite2013,ahmedNumericalStudiesHigher2017} for the nonstationary case. The efficient solution of the arising algebraic systems of equations with a huge number of unknowns is no less difficult. Often, the linear solver represents the limiting factor for the level of mesh resolution and its number of degrees of freedom. If higher order time discretizations are used, this even puts an additional facet of complexity on the structure of the resulting linear systems and their numerical solution. Here, a geometric multigrid (GMG) method that is used as a preconditioner for generalized minimal residual (GMRES) iterations is proposed and analyzed computationally for higher order space-time finite element discretizations of the Navier--Stokes equations in two and three space dimensions. Its parallel implementation in the \emph{deal.II} library \cite{arndtDealIILibrary2020} is addressed.
\textblue{The application of multigrid techniques to space-time finite element approximations of the Navier—Stokes equations and the advanced implementational issues of these algorithms yields the innovation of this work. To the best of our knowledge, the multigrid efficiency for such systems has not been studied sufficiently yet and deserves further elucidation.}

Discretizing the Navier--Stokes system by inf-sup stable pairs of finite elements that are used in this work and applying Newton's method for the linearization leads to linear systems of equations with saddle point structure. Higher order variational time discretizations, that we implement as time marching schemes by the choice of a discontinuous in time test basis, lead to linear block systems within which each block is a saddle point problem itself (cf.~\eqref{eq:newton_system_matrix_dG_1}). Solving indefinite saddle point problems has been studied intensively in the literature, cf.~\cite{benziNumericalSolutionSaddle2005,elmanFiniteElementsFast2014}. Using a direct solver is a suitable approach for problems of small dimensions. Due to the increase in computational costs and memory, two-dimensional problems that are of interest in practice or even three-dimensional simulations are not feasible for direct linear solvers, even if parallelism is used. In such cases, Krylov subspace methods \cite{saadGMRESGeneralizedMinimal1986} or multigrid schemes \cite{turekEfficientSolversIncompressible1999} are typically applied. A classical choice for a Krylov method is the (flexible) generalized minimal residual (GMRES) method. One drawback of the GMRES solver is that an additional amount of memory is allocated in each iteration. As a remedy, a restart that typically leads to a lower rate of convergence can be used; cf.~\cite{wathenPreconditioning2015}. To improve the convergence of the GMRES method, a preconditioner is typically applied within the GMRES iterations. If the density and viscosity of the flow are constant, the "pressure convection-diffusion" (PCD) block preconditioner results in a mesh independent convergence behavior; cf.~\cite{kayPreconditionerSteadyStateNavier2002,elmanFiniteElementsFast2014,cyrStabilizationScalableBlock2012}). This holds  
at least for flow problems with low to medium Reynolds numbers; \cite{olshanskiiPressureSchurComplement2007}.
As an alternative, algebraic multigrid (AMG) methods can be used for preconditioning the GMRES method. For saddle point problems, the AMG preconditioner can not be applied in a natural way. In \cite{notayNewAlgebraicMultigrid2016}, its application becomes  feasible by an appropriate transformation of the underlying saddle point problem.

In challenging benchmarks, e.g., for flow around a cylinder \cite{schaferBenchmarkComputationsLaminar1996}), geometric multigrid (GMG) methods have proven to belong to the most efficient solvers, that are currently available, cf.~\cite{johnNumericalPerformanceSmoothers2000}.
Numerical studies showed that the performance, robustness and efficiency of the GMG methods can be improved further, if they are applied as preconditioners for Krylov iterations; cf.~\cite{johnHigherOrderFinite2002}. This combination has shown to work robustly even in challenging three-dimensional simulations \cite{johnEfficiencyLinearizationSchemes2006}.
GMG methods have also been applied successfully in two space dimensions along with higher order space-time finite element discretizations of convection-diffusion equations \cite{vandervegtHpMultigridSmootherAlgorithm2012,vandervegtHpMultigridSmootherAlgorithm2012a} or the Navier--Stokes equations \cite{hussainEfficientNewtonmultigridSolution2014,hussainNoteAccurateEfficient2012,hussainEfficientStableFinite2013}. Applying the GMG method involves several complexities. Firstly, one needs to store the problem structure on various mesh levels and transfer information from finer mesh levels to coarser mesh levels by restriction operators and vice versa by prolongation operators. Secondly, parallel assembly routines add a further layer of complexity to the previous ones.
Finally, the key ingredient of the GMG method is the choice of the smoother, which damps high frequency errors on successively coarser mesh levels. The classical Gauss-Seidel smoother is not applicable to the Navier--Stokes equations, due to the saddle point structure of the discrete system. Two popular choices for this kind of problems are the Vanka type smoothers \cite{vankaBlockimplicitMultigridCalculation1986} and the  Braess–Sarazin type smoothers \cite{braessEfficientSmootherStokes1997}. In numerical studies, Vanka type smoothers have shown to outperform the Braess-Sarazin ones \cite{johnNumericalPerformanceSmoothers2000}.

In this work we propose a GMG approach based on a local (cell-based) Vanka smoother for higher order discontinuous Galerkin approximations in time. An inf-sup stable pair of finite element spaces with discontinuous pressure approximation is used for the discretization in space. This GMG method is built in the state of the art, multi-purpose finite element toolbox {deal.II} (cf.~\cite{arndtDealIILibrary2020} and \cite{clevengerFlexibleParallelAdaptive2021,kronbichlerFastMassivelyParallel2018}) along with the linear algebra package Trilinos \cite{thetrilinosprojectteamTrilinosProjectWebsite2020}. Efficient data structures are provided.
The deal.II library is enhanced in such a way that a parallel assembly and application of a cell-based Vanka smoother become feasible.
We note that some geometric multigrid support, that was used in \cite{kanschatMultigridMethodsTextdiv2015} to implement a GMG based preconditioner using $H^{\text{div}}$-conforming elements for the Stokes problem, was already provided in deal.II. The robustness and efficiency of our GMG method in terms of a grid-independent convergence of the preconditioned GMRES iterations is studied computationally for the popular benchmark problems of flow around a cylinder in two and three space dimensions.
For this, we explicitly note that parallel multigrid iterations for challenging three-dimensional flow problems do by far not meet a standard nowadays. This is underlined by the fact that the three-dimensional benchmark problem of flow around a cylinder \cite{schaferBenchmarkComputationsLaminar1996} continues to be an open one.
Confirmed numbers for the goal quantities of the simulation are not available yet.

This work is organized as follows. In Sec.~\ref{sec:Mathematical_problem}, the prototype model problem and the notation are introduced. In Sec.~\ref{Sec:STFEM}, our space-time finite element approach for simulating the Navier--Stokes system, as well as the structure of the resulting underlying system matrix, are presented. In Sec.~\ref{sec:multigrid}, we briefly recall the geometric multigrid algorithm as well as the local Vanka smoother used in this work.
We address practical aspects of the parallel implementation of the algorithm by using the deal.II library and the linear algebra package Trilinos. In Sec.~\ref{sec:numerical_example}, we present numerical results and measure the performance properties of our proposed algorithms for the 2d and 3d DFG benchmark of flow around a cylinder \cite{schaferBenchmarkComputationsLaminar1996}.

\section{Mathematical problem, notation and discretization}
\label{sec:Mathematical_problem}

\subsection{Model problem}

\begin{figure}[h!tb]
	\centering
	\includegraphics[width=0.6\linewidth]{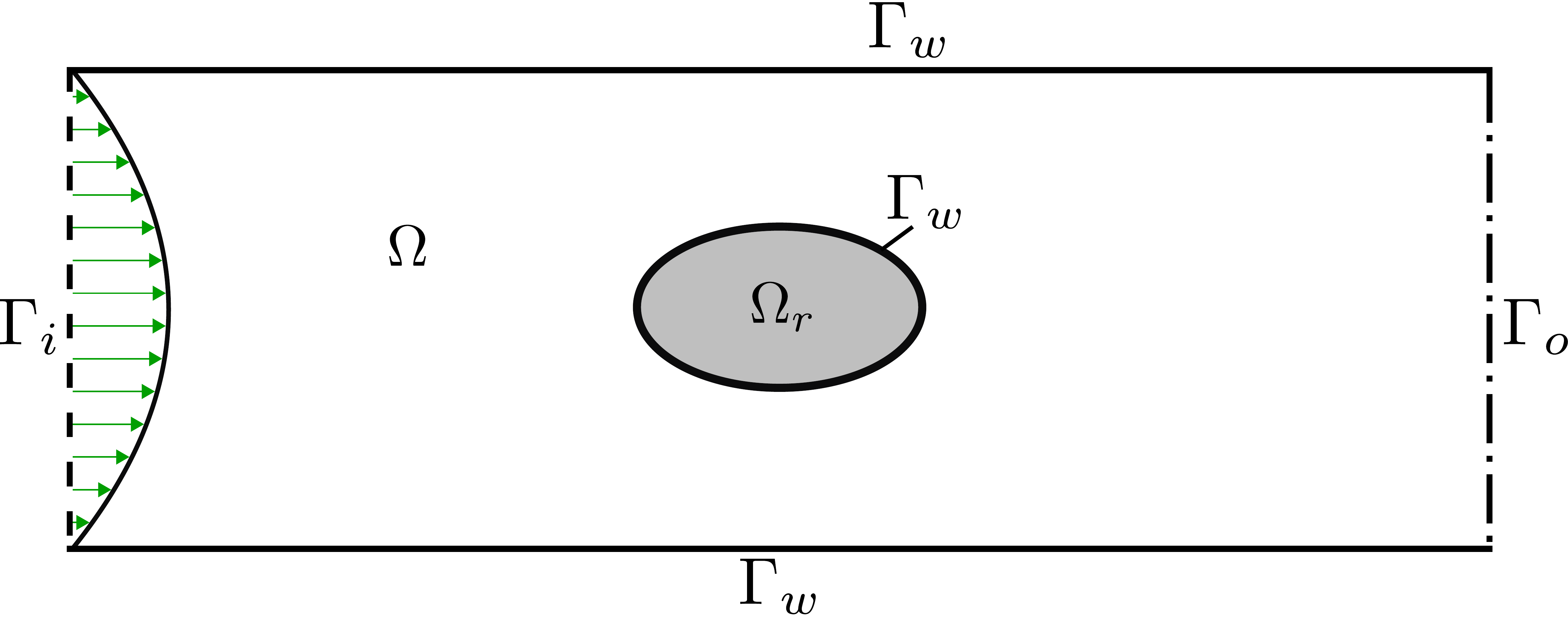}
	\caption{Problem setting for two space dimensions and corresponding notation.}%
	\label{fig:problem_overview}
\end{figure}

Without loss of generality, we consider the prototype model problem of incompressible viscous flow around a cylinder in a rectangular two- or three-dimensional domain. The two-dimensional problem configuration along with the notation of the geometrical setting is sketched in Fig.~\ref{fig:problem_overview}. We evaluate the performance properties of the proposed GMG solver for this benchmark problem (cf.\cite{schaferBenchmarkComputationsLaminar1996}).
We consider solving the Navier--Stokes equations
\begin{subequations}
	\label{eq:navier_stokes}
	\begin{alignat}{4}
		\label{eq:navier_stokes_0}
		\partial_t \vec{v} + (\vec{v} \cdot \nabla) \vec{v} - \nu \Delta \vec{v} + \nabla p	&= \vec{f}
		& \hspace*{2ex} & \text{in } \Omega \times I\,,
		\\
		\label{eq:navier_stokes_1}
		\nabla \cdot \vec{v} &= 0
		& & \text{in } \Omega  \times I\,,
		\\
		\label{eq:navier_stokes_2}
		\vec{v} & = \vec{g}
		& & \text{on } \Gamma_D \times I\,, \quad
		\\
		\label{eq:navier_stokes_5}
		\nu \nabla \vec{v} \cdot \vec{n} - \vec{n}p &= \vec 0
		& & \text{on } \Gamma_{o}  \times I\,, \quad
		\\
		\label{eq:navier_stokes_6}
		\vec{v}(0) &= \vec{v}_0 
		& & \text{in } \Omega \,.
	\end{alignat}
\end{subequations}
In \eqref{eq:navier_stokes}, $\Omega \subset \R^d$, with $d=2$ or $d=3$, is the open domain filled with fluid. We put $I=(0,T]$ for some final time $T>0$. The velocity field $\vec v$ and the pressure $p$ are the unknown variables. In \eqref{eq:navier_stokes_0}, the parameter $\nu>0$ denotes the fluid's viscosity and the right-hand side function $\vec{f}$ is a given external force.
The union of the Dirichlet boundary segments is denoted by $\Gamma_D$, such that $\Gamma_D:=\Gamma_i\cup \Gamma_w$. On $\Gamma_D$ we prescribe  the fluid velocity by a function $\vec{g}$, that prescribes an inflow profile on $\Gamma_i$ and a no slip condition on $\Gamma_w$. $\Gamma_o$ represents an outflow boundary that is modeled by the do-nothing boundary condition \eqref{eq:navier_stokes_5}; cf.~\cite{heywoodArtificialBoundariesFlux1996}.
In \eqref{eq:navier_stokes_5}, the field $\vec n$ is the outer unit normal vector. In \eqref{eq:navier_stokes_6}, the function $\vec{v}_0$ denotes the prescribed initial velocity. In our numerical experiments presented in Sec.~\ref{sec:numerical_example}, the (time-independent) rigid domain $\Omega_r$ is a sphere in two space-dimensions or a cylinder in three space-dimensions.

We assume that a sufficiently regular solution of the system \eqref{eq:navier_stokes} exists such that higher order approaches become feasible. For the existence, uniqueness and regularity of solutions to the Navier--Stokes system, including the regularity for $t\rightarrow 0$ under realistic assumptions about the data, we refer to the broad literature in this field; cf.\ \cite{johnFiniteElementMethods2016,heywoodFiniteElementApproximation1982} and the references therein. 

\subsection{Notation}
\label{sec:Problem_Notation}

Here, we introduce the function spaces that are used in this work to present our space-time approach for the Navier--Stokes system \eqref{eq:navier_stokes}. By $L^2(\Omega)$ we denote the function space of square integrable functions on the fluid domain $\Omega$ while $H^1(\Omega)$ is the usual Sobolev space of functions in $L^2(\Omega)$ which have first order weak derivatives in $L^2(\Omega)$. Further, $\langle \cdot, \cdot \rangle$ is the standard inner product of $L^2(\Omega)$. We define the subspace of $L^2(\Omega)$ with mean zero $L^2_0(\Omega):= \{v\in L^2(\Omega) \mid \int_\Omega v \d x = 0\}$ and the subspace of $H^1(\Omega)$ of functions with zero boundary values (in the sense of traces) on the portion $\Gamma_D\subset \partial \Omega$ of the boundary $\partial \Omega$ of $\Omega$ as $H^1_{0,\Gamma_D}(\Omega)$. Its dual space is denoted by $H^{-1}(\Omega)$. Finally, by $H^{1/2}(\Gamma_D)$ we denote the space of all traces on $\Gamma\subset \partial \Omega$ of functions in $H^1(\Omega)$. For vector-valued functions we write those spaces bold. 
%
%


\subsection{Space discretization}
\label{Subsec:SpaceDisc}

Let $\mathcal{T}_h=\{T\}$ be a family of shape-regular decompositions of the $\Omega$ (cf.\ Fig.~\ref{fig:problem_overview}) into (open) quadrilaterals $K$ with maximum cell size $h$.
For $r \geq 0$, let $\hat{\mathbb{Q}}_r$ denote the space of polynomials of degree at most $r$ in each variable and $\hat{\mathbb{P}}_r$ the space of at most degree $r$. We put 
\begin{equation*}
	\mathbb{Q}_r(T) \coloneqq \left\{
	v_h{}_|{}_T = \hat{v}_h \circ F_T^{-1} : \hat{v}_h \in \hat{\mathbb{Q}}_r
	\right\} \,,
	\quad
	\mathbb{P}_r(T) \coloneqq \left\{
	v_h{}_|{}_T = \hat{v}_h \circ F_T^{-1} : \hat{v}_h \in \hat{\mathbb{P}}_r
	\right\} \,,
\end{equation*}
with the reference mapping $F_T$ from the reference cell $\widehat T$ to element $T$.
In our computations presented in Sec.~\ref{sec:numerical_example}, only affine linear mappings are used.
We define the finite element spaces 
\begin{align*}
	\label{Eq:DefHh}
	H_h^{r}&\coloneqq \left\{v_h \in C({\overline{\Omega}}) \mid v_h{}_|{}_T\in \mathbb{Q}_r(T) \; \forall \; T \in \mathcal{T}_h \right\}\,, \\
	H_{h, \text{disc}}^{r}&\coloneqq \left\{v_h \in L^2({\Omega}) \mid v_h{}_|{}_T\in \mathbb{P}_r(T) \; \forall \; T \in \mathcal{T}_h \right\}\,.
\end{align*}
For the construction and effective application of the cell-based Vanka smoother, a local velocity-pressure coupling is required. To establish such a coupling, a discontinuous in space pressure approximation is applied.
For the spatial approximation of the velocity and pressure variable we use the conforming, inf-sup stable finite element pair that is given by
\begin{equation}
	\label{Def:Vh0}
	\vec V_h = H_h^{r} \times H_h^{r} \,, \quad Q_h = H_{h,\text{disc}}^{r-1}
\end{equation}
for some natural number $r\geq 2$, cf.~\cite{johnFiniteElementMethods2016,matthiesInfsupConditionMapped2002}. All the numerical experiments that are presented in Sec.~\ref{sec:numerical_example} were done by the choice \eqref{Def:Vh0} of the discrete function spaces.

The space of weakly divergence free functions is denoted by
\begin{equation*}
	\vec V_h^{\text{div}} = \left\{ \vec v_h \in \vec V_h \mid \langle \nabla \cdot \vec 
	v_h,q_h\rangle  = 0 \; \text{for all } q_h \in Q_h \right\} \,.
\end{equation*}
Finally, we define the spaces  
\begin{equation*}
	\label{Eq:DefSolSpa1_h}
	V_{I,h} := \{\vec v_h \in L^2(I;\vec V_h) \mid \partial_t \vec v \in L^2(I;\vec V_h)\}\,, \qquad L^2_{0,I,h} : = L^2(I;Q_h)\,.	
\end{equation*}

\textblue{
For the treatment of Dirichlet boundary conditions by Nitsche's method (cf.~\cite{beckerMeshAdaptationDirichlet2002}) we introduce the bilinearform  $B_{\Gamma_D} : \vec H^{1/2}(\Gamma_D) \times (\vec V_h \times Q_h) \rightarrow \R$ by 
\begin{equation}
	\label{Def:B_GamD}
	\begin{aligned}
		B_{\Gamma_D}(\vec w,(\vec \psi_h,\xi_h)) : = &  - \langle \vec w, \nu \nabla \vec \psi_h \cdot \vec n + \xi_h  \vec n  \rangle_{\Gamma_D}
		\\[1ex]
		& \quad + \gamma_1 \nu \langle h^{-1} \vec w , \vec 
		\psi_h  \rangle_{\Gamma_D}  + \gamma_2 \langle h^{-1} \vec w \cdot \vec n, \vec \psi_h \cdot \vec n 
		\rangle_{\Gamma_D} 
	\end{aligned}
\end{equation}
for $\vec w \in \vec H^{1/2}(\Gamma_D)$ and $(\vec \psi_h,\xi_h) \in \vec V_h \times Q_h$, where $\gamma_1>0$ and $\gamma_2> 0$ are numerical (tuning) parameters for the penalization. In  \cite{agerNitschebasedCutFinite2019,winterNitscheCutFinite2018}, their choice in the range of $(10, 100)$ is recommended. In our simulations presented in Sec.~\ref{sec:numerical_example}, we put $\gamma_1 = \gamma_2 = 35$. 
Next we define the semilinear form $A_h: (\vec V_h\times Q_h)\times (\vec V_h\times Q_h)\rightarrow \R$ by
\begin{equation}
\label{Def:Ah}
\begin{aligned}
	A_h((\vec v_h,p_h),(\vec \psi_h,\xi_h)) \coloneqq &
	\langle (\vec v_h \cdot \vec \nabla) \vec v_h, \vec \psi_h \rangle
	+ \nu \langle \nabla \vec v_h , \nabla \vec \psi_h  \rangle
	-\langle p_h, \nabla \cdot \vec \psi_h \rangle
	+ \langle \vec \nabla \cdot \vec v_h, \xi_h \rangle \\
	& -
	\langle \nu \nabla \vec 
	v_h \cdot \vec n - p_h \vec n, \vec \psi_h\rangle_{\Gamma_D} + B_{\Gamma_D}(\vec v_h,\vec 
	\phi_h)
\end{aligned}
\end{equation}
for $(\vec v_h,p_h)\in \vec V_h\times Q_h$ and $(\vec \psi_h,\xi_h) \in \vec V_h\times Q_h$. The linear form $L_h: (\vec V_h\times Q_h) \rightarrow \R$ is defined by
\begin{equation}
	\label{Def:Lh}
	L_h((\vec \psi_h,\xi_h);\vec f, \vec g)	:= L(\vec \psi_h; \vec f) + B_{\Gamma_D}(\vec g,(\vec \psi_h,\xi_h)) 
\end{equation}
for $(\vec \psi_h,\xi_h) \in \vec V_h\times Q_h$.
}

\subsection{Time discretization}

For the time discretization, we decompose the time interval $I=(0,T]$ into $N$ subintervals $I_n=(t_{n-1},t_n]$, $n=1,\ldots,N$, where $0=t_0<t_1< \cdots < t_{N-1} < t_N = T$ such that $I =\bigcup_{n=1}^N I_n$ and $I_n \cap I_m = \emptyset$ for $n \neq m,\,m,n = 1,\ldots,N$.
We put $\tau = \max_{n=1,\ldots, N} \tau_n$ with $\tau_n = t_n-t_{n-1}$. Further, the set $\mathcal{M}_\tau := \{I_1,\ldots, I_N\}$ of time intervals is called the time mesh. For a Banach space $B$ of functions defined on the time-independent domain $\Omega$ and any $k\in \N_0$, we let 
\begin{align}
	\label{Def:PkIB}
	\P_k(I_n;B)
	= \Big\{w_\tau : I_n \to B \mid w_\tau(t) = \mbox{$\sum\limits_{j=0}^k$}
	W_j \, t^j \; \forall t \in I_n\,, \;
	W_j \in B\; \forall j \Big\}\,.
\end{align}
For an integer $k\in \N_0$, we put 
\begin{equation}
	\label{Def:Xtauk}
	X_\tau^{k} (B)
	\coloneqq
	\left\{w_\tau \in L^2(I;B) \mid w_\tau|_{I_n} \in
	\P_{k}(I_n;B)\; \forall I_n\in \mathcal{M}_\tau\,, \; w_\tau (0) \in B \right\}\,. 
\end{equation}

\section{Space-time finite element discretization}
\label{Sec:STFEM}

\textblue{
For the discretization of the Navier--Stokes system \eqref{eq:navier_stokes} by space-time finite element methods, we follow the lines of \cite{turekEfficientSolversIncompressible1999,johnNumericalPerformanceSmoothers2000,anselmannHigherOrderGalerkin2020,anselmannCutFiniteElement2022} and consider the time marching scheme, that consists of solving the sequence of the following local problems defined on the subinterval $I_n$:
\begin{problem}
	\label{Prob_FullyDisc}
	Let $\vec f \in  L^2(I;\vec H^{-1}(\Omega))$ and $\vec v_{0,h} \in \vec V^{\text{div}}_h$ be given. For $n=1,\ldots, N$, and given $\vec v_{\tau,h}{}_{|I_{n-1}} \in \mathbb P_k(I_{n-1}; \vec V_h )$ for $n>1$ and $\vec v_{\tau,h}{}_{|I_{n-1}}(t_{n-1}^-) := \vec v_{0,h}$ for $n=1$, find $(\vec v_{\tau,h},p_{\tau,h})\in \mathbb P_k(I_n; \vec V_{h}) \times \mathbb P_k(I_n;Q_h)$, such that 
	\begin{multline}
		\label{Prob_FullyDisc_1}
		\int_{t_{n-1}}^{t_{n}} \langle\partial_t \vec v_{\tau,h}, \vec \psi_{\tau,h} \rangle +
		A_h((\vec v_{\tau,h},p_{\tau,h}),(\vec \psi_{\tau,h},\xi_{\tau,h})) \d t + \langle \vec v_{\tau,h}(t_{n-1}^+), \vec \psi _{\tau,h}(t_{n-1}^+) \rangle \\ = \int_{t_{n-1}}^{t_{n}}
		L_h(\psi_{\tau,h}; \vec f, \vec g) \d t + \langle \vec v_{\tau,h}(t_{n-1}^-), \vec \psi _{\tau,h}(t_{n-1}^+) \rangle
	\end{multline}
	for all $(\vec \psi_{\tau,h},\xi_{\tau,h}) \in \mathbb P_k(I;\vec V_h)\times \mathbb P_k(I_n;Q_h)$.
\end{problem}
This discretization features a discontinuous Galerkin in time method, with piecewise polynomials of order $k\in \N_0$.
Using a discontinuous method in the time domain has the advantage, that no initial value for the pressure is needed. Continuous in time methods require a discrete initial pressure value for the unique definition of its full trajectory.
Such an initial value, that guarantees the optimal order of convergence of the velocity and pressure variables for all $t\in I$, is not available.
A remedy is the application of extrapolation techniques; cf.~\cite{hussainNoteAccurateEfficient2012}.
Further, discontinuous Galerkin methods offer stronger stability properties since they are known to be strongly $A$-stable.
}

For the construction of the GMG method in Subsec.~\ref{sec:cell_based_vanka} we discuss the algebraic counterpart of Eq.~\eqref{Prob_FullyDisc_1} more thoroughly. Firstly, we represent the unknown discrete functions $(\vec v_{\tau,h},p_{\tau,h})\in L^2(I_n; \vec V_{h}) \times L^2(I_n;Q_h)$ in a temporal basis $\{\chi_l\}_{l=0}^k$ of $\mathbb P_k (I_n;\R)$ by means of  
\begin{align}
	\label{eq:time_discrete_ansatz}
	{v}_{\tau,h,i|I_n}(\mat x, t)
	&=
	\sum_{l=0}^{k} {v}^{n,l}_i
	(\mat{x})\chi_{n,l}(t)\,,\;\; \text{for } i\in \{1,\ldots, d\}\,,
	&
	p_{\tau,h|I_n}(\mat x, t)
	&=
	\sum_{l=0}^{k}
	p^{n,l}
	(\mat{x})\chi_{n,l}(t)\,,
\end{align}
with coefficient functions $\vec{v}^{n,l} = (v^{n,l}_1,\ldots, v^{n,l}_d)^\top\in \vec V_h$ and $p^{l,n} \in Q_h$, where $\vec v_{\tau,h} = (\vec v_{\tau,h,1},\ldots,$ $\vec v_{\tau,h,d})^\top$ for $t\in I_n$. For the basis $\{\chi_l\}_{l=0}^k$ we choose the Lagrange interpolants with respect to the $k+1$ Gauss--Radau quadrature nodes of $I_n$. Appreciable of the Gauss--Radau quadrature formula is that the end point of the subinterval $I_n$ is a quadrature node, which simplifies the evaluation of the second term on the right-hand side of \eqref{Prob_FullyDisc_1}. Letting 
\begin{equation}
\label{eq:spacebasis}
H_h^r = \operatorname{span}\{\psi_1,\ldots, \psi_R\} \quad \text{and}\quad 
Q_h = \operatorname{span}\{\xi_1,\ldots , \xi_S\}\,,
\end{equation}
the coefficient functions $\vec{v}^{n,l}\in \vec V_h$ and $p^{n,l}\in \vec Q_h$ of \eqref{eq:time_discrete_ansatz} admit the representation 
\begin{equation*}
\label{eq:expcoeffunc}
{v}^{n,l}_i(\vec x) = \sum_{r=1}^R {v}^{n,l}_{i,r} \, \psi_r(\vec x)
\quad \text{and}\quad 
p^{n,l}(\vec x) = \sum_{s=1}^S {p}^{n,l}_s \, \xi_s(\vec x)
\end{equation*}
with the vectors of unknown coefficients 
\begin{equation}
\label{eq:coeffvec}
\vec {v}^{n,l}_{i} = ({v}^{n,l}_{i,1},\ldots,{v}^{n,l}_{i,R})^\top \in \R^R
\quad \text{and}\quad
\vec {p}^{n,l} = ({p}^{n,l}_1,\ldots,{p}^{n,l}_S)^\top \in \R^S\
\end{equation}
for all degrees of freedom in time in $I_n$ with $l=0,\ldots, k$. Clearly, the vectors $\vec {v}^{n,l}_{i}$ denote the coefficients of the velocity component functions ${v}^{n,l}_i$, with $i=1,\ldots ,d$, with respect to the spatial basis $\{\psi_r\}_{r=1}^R$. 

Defining now the vector $\vec X^n\in \R^{(k+1)\times (d\cdot R+S)}$ of unknown coefficients for the solution of \eqref{Prob_FullyDisc_1} in the subinterval $I_n$ by 
\begin{equation}
\label{eq:totdofs}
\vec X_n = (\vec {v}^{n,0}_{1},\ldots,\vec {v}^{n,0}_{d},\vec {p}^{n,0},\ldots, \vec {v}^{n,k}_{1},\ldots,\vec {v}^{n,k}_{d},\vec {p}^{n,k})^\top \in \R^{(k+1)\times (d\cdot R+S)}\,, 
\end{equation}
we recover the variational equation \eqref{Prob_FullyDisc_1} in an algebraic form as 
\begin{equation}
\label{eq:nonlineq}
\vec F_n(\vec X_n) = \vec 0 \,,
\end{equation} 
for a suitably defined nonlinear function $\vec F_n: \R^{(k+1)\times (d\cdot R+S)} \rightarrow \R^{(k+1)\times (d\cdot R+S)}$. We refer to \cite{anselmannHigherOrderGalerkin2020} for the explicit derivation of the algebraic system of a related space-time finite element approximation of the Navier--Stokes system. To solve the nonlinear problem \eqref{eq:nonlineq}, we use Newton's method such that the linear system 
\begin{equation}
\label{eq:linNewEq}
\vec J_{n}^m \vec D_n^m = \vec Q_n^m\,, 
\end{equation} 
with the Jacobian matrix and right-hand side vector 
\begin{equation}
\label{eq:definition_jacobian_global}
\vec J_{n}^m  = \dfrac{\partial \vec F_{n}}{\partial \vec X_n}(\vec X_n^m)\quad \text{and} \quad \vec Q_n^m = \vec F_n(\vec X_n^m)\,,
\end{equation}
has to be solved in each Newton iteration $m$ for the new iterate
\begin{equation}
	\label{eq:linNewEqU}
	\vec X_n^{m+1} = \vec X_n^m + \vec D_n^m\,.
\end{equation} 
For brevity, an explicit form of the Jacobian matrix $\vec J_{n}^m $ is not given here. For the sake of clarity, we restrict ourselves to presenting the block structure of $\vec J_{n}^m$ for the polynomial order in time $k = 1$ only. This corresponds to the dG(1) method for the time discretization. For this, we get that
\begin{equation}
	\label{eq:newton_system_matrix_dG_1}
	\vec J_{n}^m=
	\begin{pmatrix}
		\mat{F}_1 & \mat{B}_1^\top & \mat{F}_2 & \mat{B}_2^\top \\[1ex]
		-\mat{B}_1 & \mat{0} & -\mat{B}_2 & \mat{0} \\[1ex]
		\mat{F}_3 & \mat{B}_3^\top & \mat{F}_4 & \mat{B}_4^\top \\[1ex]
		-\mat{B}_3 & \mat{0} & -\mat{B}_4 & \mat{0}
	\end{pmatrix} .
\end{equation}
In \eqref{eq:newton_system_matrix_dG_1}, the partitioning of the vector of unknowns $\vec D_n^m$ of the corresponding system  \eqref{eq:linNewEq} is then given by
\begin{equation*}
	\vec D_n^m = (\vec {d_{\vec v}}^{n,0},\vec {d_p}^{n,0},\vec {d_{\vec v}}^{n,1},\vec {d_p}^{n,1})^\top \in \R^{2(d \cdot R + S)}\,,
\end{equation*}
where $\vec {d_{\vec v}}^{n,l}$ and $\vec {d_{p}}^{n,l}$, with $l\in \{0,1\}$, denote the components of $\vec D_n^m$ related to velocity and pressure unknowns, respectively. 
For a more detailed derivation of the system matrix of  \eqref{eq:newton_system_matrix_dG_1} and the definition of the submatrices $\vec B_i$ and $\vec F_i$ we refer to \cite{anselmannHigherOrderGalerkin2020} again.
\textblue{
\Cref{tab:used_indices} summarizes the indices, used throughout this section.
\begin{table}[!ht]
	\caption{\textblue{Summary of the indices used to describe the fully discrete problem.}}
	\centering
	\textblue{
	\begin{tabular}{ c l l}
		\toprule
		index 	& {range} 				& {explanation} \\ \hline
		$i$ 	& $[1,\ldots,d]$ 		& spatial velocity components, $d$ = spatial dimension \\
		$n$ 	& $[1,\ldots,N]$ 		& time interval $I_n$\\
		$l$ 	& $[0,\ldots,k]$		& local DoFs in time on $I_n$ \\
		$r$ 	& $[1,\ldots,R]$		& spatial DoFs for $\vec{v}$ \\
		$s$ 	& $[1,\ldots,S]$		& spatial DoFs for $p$ \\
		$m$  	& $\geq 1$				& counter for Newton iteration \\
		\bottomrule
	\end{tabular}
	}
	\label{tab:used_indices}
\end{table}
}

To enhance the range of convergence of Newton's method, a damped version using an additional linesearch technique is applied for solving \eqref{eq:nonlineq}. Alternatively, a "dogleg approach" (cf., e.g. \cite{pawlowskiInexactNewtonDogleg2008}), that belongs to the class of trust-region methods and offers the advantage that also the search direction, not just its length, can be adapted to the nonlinear solution process, was implemented and tested. Both schemes require the computation of the Jacobian matrix of the algebraic counterpart of Eq.~\eqref{Prob_FullyDisc_1}. In the dogleg method multiple matrix-vector products with the Jacobian matrix have to be computed. Since the Jacobian matrix is stored as a sparse matrix, the products can be computed at low computational costs.
From the point of view of convergence, both methods yield a superlinear convergence behavior. In our numerical examples of Sec.~\ref{sec:numerical_example}, both modifications of Newton's method lead to comparable results. In our computational studies, we did not observe any convergence problems for these nonlinear solvers. To solve the linear systems \eqref{eq:linNewEq} of the Newton iteration, we use a flexible GMRES Krylov subspace method \cite{saadGMRESGeneralizedMinimal1986} with a GMG preconditioner based on a local Vanka smoother \cite{turekEfficientSolversIncompressible1999}.
The GMG approach is presented in the next section.


\section{A parallel geometric multigrid preconditioner}
\label{sec:multigrid}

During the last decades numerous methods for solving the algebraic linear systems resulting from the discretization of the Navier--Stokes equations have been developed and studied. GMG methods seem to be among the best classes of solvers that are currently available; cf., e.g., \cite{johnNumericalPerformanceSmoothers2000}. Space-time finite element methods have recently attracted researchers' interest strongly. Their application puts an additional complexity to the solution of the linear systems due to their more complex block structure; cf.\ e.g.,  \cite{hussainEfficientNewtonmultigridSolution2014,anselmannHigherOrderGalerkin2020,bastingPreconditionersDiscontinuousGalerkin2017}.
Here, we use the GMG method as a preconditioner for Krylov subspace iterations, which is a standard concept for the efficient solution of high-dimensional linear systems arising from the discretization of partial differential equations. The core of GMG methods is the smoother. We propose a cell-based Vanka smoother that is adapted to the space-time finite element approach.

Even though the basic concepts of GMG methods have become standard, their efficient implementation continues to be a challenging task. In particular, this holds if the computational power of modern parallel computer architectures has to be fully exploited. In this case, the definition of data structures and the memory usage become of utmost importance. Moreover, trends like adaptive space-time finite element methods (cf.~\cite{kocherEfficientScalableData2019}) further complicate their implementation. These issues induce an ongoing research about GMG methods; cf., e.g., \cite{clevengerFlexibleParallelAdaptive2021}. For our simulations we use the \emph{deal.II} finite element toolbox \cite{arndtDealIILibrary2020}. Details of our implementation of the GMG approach in this platform are addressed in the sequel as well. The concepts are flexible enough and can be transferred to similar software tools.

\subsection{Key idea of the geometric multigrid method}

\begin{figure}[!htb]
	\centering
	\includegraphics[width=1.0\textwidth,keepaspectratio]
	{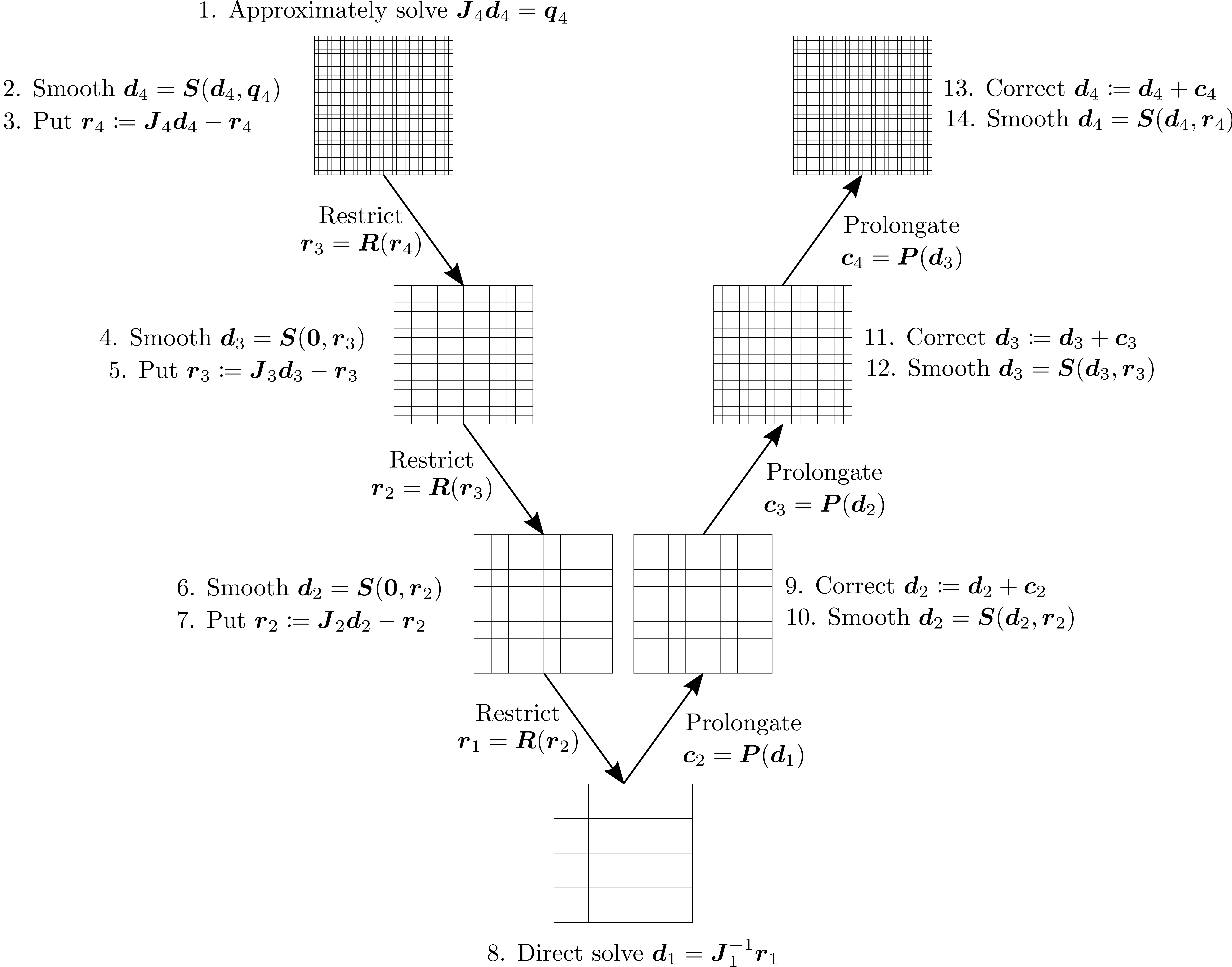}
	\caption{Structure of a single multigrid V-cycle for a hierarchy of four grid levels.} 
	\label{Fig:MG_overview}
\end{figure}

To sketch briefly the basic principles of GMG iterations and fix our notation, we consider the linear system \eqref{eq:linNewEq}, that is rewritten in the simpler, index-free form 
\begin{equation}
	\label{eq:LinSysMG}
	\vec J \vec d = \vec r\,,
\end{equation}
with the right-hand side vector $\vec r$ and the Jacobian matrix $\mat J$. The key idea of the GMG method, that is sketched in Fig.~\ref{Fig:MG_overview}, is to construct a hierarchical sequence of finite element spaces $V_h^g$, with $g=1,\ldots, G$, that are embedded into each other, such that $V_h^{1} \subset V_h^{2} \subset \ldots \subset V_h^{G}$, and correspond to different grid levels $\mathcal T_{h_g}$ with mesh sizes $h_g$, for $g=1,\ldots, G$ , of decompositions of the domain $\Omega$. Instead of solving the linear system \eqref{eq:LinSysMG} on the finest grid level $\mathcal T_{h_G}$ entirely, the idea is to smooth only high frequency errors of an initial guess to the solution $\vec d= \vec d_G$ of \eqref{eq:LinSysMG} on the finest grid level $\mathcal T_{h_g}$ with $g=G$. Clearly, on level $g=G$, the right-hand side vector $\vec r=\vec r_G$ corresponds to the right-hand side vector $\vec Q_n^k$ of the Newton system \eqref{eq:linNewEq}. Now, smoothing is done by the application of the local Vanka operator $\vec S_G$. Then, the resulting residual $\vec r_G$ of \eqref{eq:LinSysMG} for the computed approximation of $\vec d_G$ is restricted to next coarser mesh level $\mathcal T_{h_g}$, with $g=G-1$, which yields the right-hand side vector $\vec r = \vec r_{G-1}$. On level $G-1$, the high frequency errors of an initial guess (given by the null vector $\vec 0$) to the solution's correction $\vec d_{G-1}$ on $\mathcal T_{h_g}$, with $g=G-1$, is smoothened by the application of the local Vanka operator again. These operations of restricting the residual to the next coarser grid and smoothing on this level the error in the solution of the defect equation is recursively repeated until the coarsest mesh level $\mathcal T_{h_g}$, with $g=1$, is reached. On this level, typically a direct solver is used to compute the corresponding defect correction $\vec d_1$. Afterwards the computed defect correction of the coarsest level is prolongated to the next finer grid level $\mathcal T_{h_g}$, with $g=2$, and used to update the defect correction $\vec d_2$. On this level, the defect correction $\vec d_2$ is then smoothed again and, finally, prolongated to the next coarser mesh level $\mathcal T_{h_g}$, with $g=3$. These operations of prolongating successively the coarse grid correction and smoothing the modified defect correction are continued until the finest grid level $\mathcal T_{h_g}$, with $g=G$, is reached, where after the final smoothing an updated solution is obtained. This GMG approach is summarized in Algorithm~\ref{alg:multigrid_general}.

\begin{algorithm}[H]
	\SetAlgoLined
	\If{$g == 1$}{
		$\vec{d}_1 = \mat{J}_1^{-1} \vec{r}_1$ \tcp*[r]{Direct coarse solver}
	}
	\Else{
		$\vec{d}_g = \vec{S}(\vec{d}_g, \vec{r}_g)$ \tcp*[r]{Pre-smooth}
		$\vec{r}_{g} \coloneqq \mat{J} \vec{d}_g - \vec{r}_g$ \tcp*[r]{Compute residuum}
		$\vec{r}_{g-1} = \vec{R}(\vec{r}_g)$ \tcp*[r]{Restrict residuum}
		Multigrid($\vec{0}, \vec{r}_{g-1}, g-1$) \tcp*[r]{Recursively call this function}
		$\vec{c}_g = \vec{P}(\vec{d}_g)$ \tcp*[r]{Prolongate correction}
		$\vec{d}_g \coloneqq \vec{d}_g + \vec{c}_g$ \tcp*[r]{Correct solution}
		$\vec{d}_g = \vec{S}(\vec{d}_g, \vec{r}_g)$ \tcp*[r]{Post-smooth}
	}
	\caption{Recursive algorithm Multigrid($\vec{d}_g, \vec{r}_g, g$)}
	\label{alg:multigrid_general}
\end{algorithm}

For our implementation of the GMG approach and the simulations presented in Sec.~\ref{sec:numerical_example}, we use the \emph{deal.II} finite element toolbox \cite{arndtDealIILibrary2020} along with the direct, parallel \emph{SuperLU\_Dist} solver \cite{liSuperLUDISTScalable2003}. Our code is based on the contributions of \cite{clevengerFlexibleParallelAdaptive2021} to this open source framework and expands their work by a parallel, cell-based Vanka smoother. For the restriction and prolongation steps in parallel computations, the deal.II classes \emph{MultiGrid} and \emph{MGTransferPrebuilt} are used.
The latter implements the prolongation between grids by interpolation and applies the transpose operator for the restriction. The core of our GMG approach is the smoother. This operator has to be efficient in smoothing high frequency errors. Further, since the smoother is applied frequently (cf.\ Fig.~\ref{Fig:MG_overview}), this demands for its performant and scalable implementation, utilizing multiple processors, such that the hardware's potential is fully exploited. Our implementation of the smoother is presented more in detail below.


\subsection{A parallel, cell-based Vanka smoother}
\label{sec:cell_based_vanka}

\begin{figure}[!htb]
	\centering
	\includegraphics[width=0.2\textwidth,keepaspectratio]
	{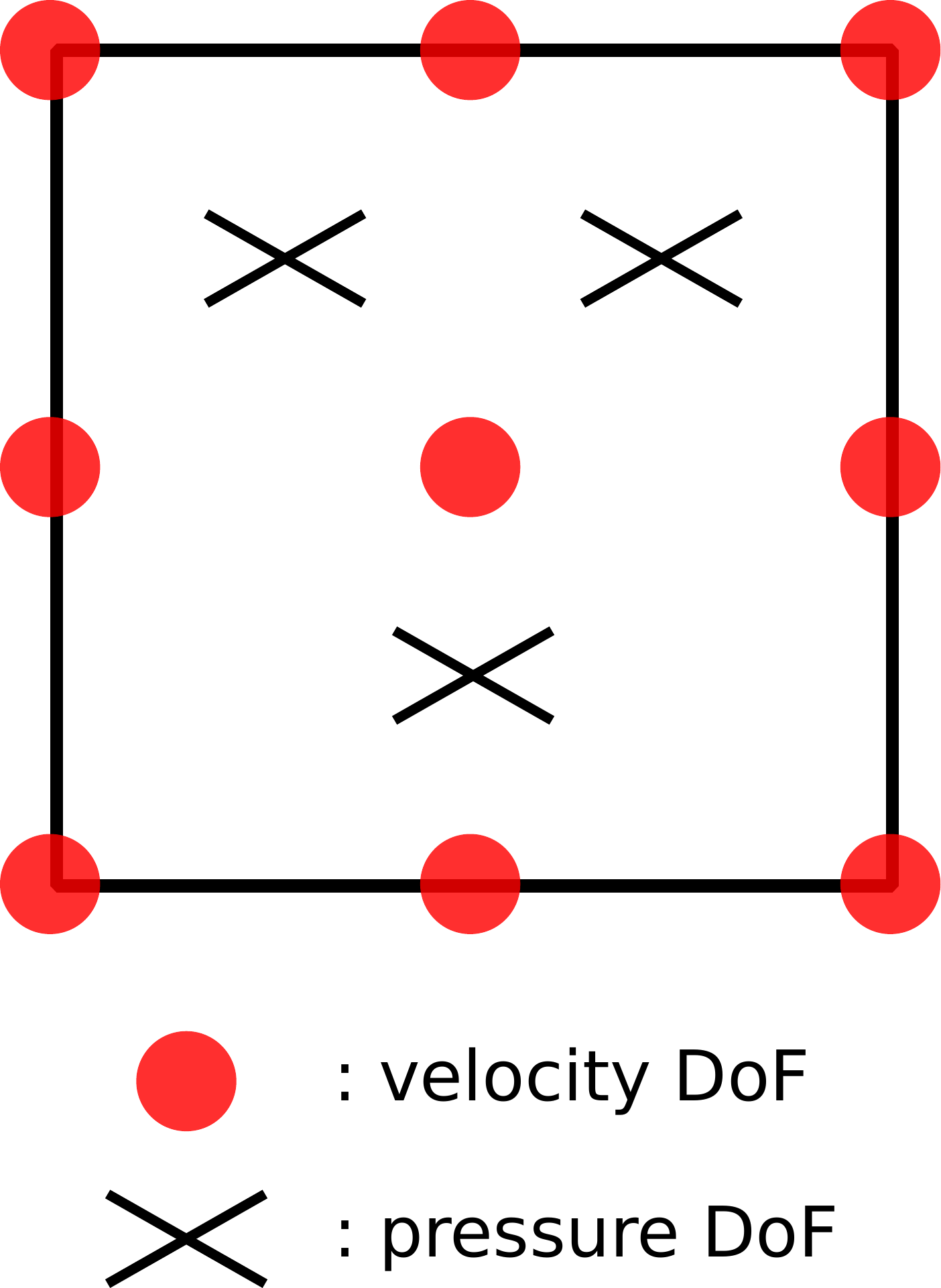}
	\caption{Degrees of freedom for the continuous/discontinuous  $\mathbb{Q}_2$--$\mathbb{P}_{1}^{\text{disc}}$ pair of elements, see \cite{matthiesInfsupConditionMapped2002}.}
	\label{Fig:FE_DGP}
\end{figure}

The Newton linearized system \eqref{eq:LinSysMG} of the fully discrete problem \eqref{Prob_FullyDisc_1} has a generalized saddle-point structure; cf.\ \cref{eq:newton_system_matrix_dG_1}.
The generalization comes through the application of the higher order discontinuous Galerkin time discretization with $k+1$ temporal degrees of freedom (cf.\ \eqref{eq:time_discrete_ansatz} and \eqref{eq:totdofs}) in time for the velocity and pressure variable within each subinterval $I_n$. Thereby, blocks of saddle point subsystems arise; cf.\ \cref{eq:newton_system_matrix_dG_1}. Standard smoothers, like the Gauss-Seidel or Jacobi method, that are often used in GMG methods, are not applicable to such systems; cf.\ \cite{benziNumericalSolutionSaddle2005a}. Vanka smoothers, that can be traced back to \cite{vankaBlockimplicitMultigridSolution1986}, offer the potential to to damp high frequency errors in the approximation of solutions to saddle point problems.
\textblue{In \cite{manservisiNumericalAnalysisVanka2006} Vanka-type solvers were analyzed for the steady Stokes problem, using Taylor-Hood elements.
Convergence was also proven for the Navier--Stokes flows in low Reynolds number regimes.}
In \cite{johnNumericalPerformanceSmoothers2000,molenaarTwogridAnalysisCombination1991} Vanka-type smoothers have demonstrated excellent performance properties for systems with weak velocity-pressure couplings. In this work, we adapt the principle of a cell-based, full Vanka smoother of \cite[p.\,460]{johnNumericalPerformanceSmoothers2000} and extend the definition of the Vanka smoother to our higher-order space-time finite element approximation of the Navier--Stokes system. Since the numerical results, reported for instance in \cite{turekNumericalStudiesVankaType2009}, show that a strong velocity-pressure coupling leads to a local violation of the continuity constraint, \eqref{eq:navier_stokes_1} we only use a discontinuous finite element space for the pressure variable, defined in Subsec.~\ref{Subsec:SpaceDisc}.
This results in the ability to use local test functions, that are defined on a single cell. 
Therefore, the mass conservation is fulfilled locally, cf.\ \cite{richterFluidstructureInteractionsModels2017,matthiesMassConservationFinite2007}.
Fig.~\ref{Fig:FE_DGP} illustrates the position of the underlying degrees of freedom for a pair of continuous/discontinuous finite elements for the velocity/pressure variables, corresponding to the case $r=2$ in the definitions of \eqref{Def:Vh0}. 

\begin{remark}
	In the deal.II finite element library, we use the \emph{FE\_DGP\qdist{ }} class to form a basis of $\mathbb{P}_{r}$.
	This basis is constructed using a set of polynomials of complete degree $r$ that form a Legendre basis on the unit square, i.\,e. they are $L^2$ orthogonal and normalized on the reference cell.
	Noteworthy, this element is not a Lagrangian one, so it is not defined by finding shape functions within the given function space that interpolate a particular set of points.
	Therefore, in \cref{Fig:FE_DGP} the pressure DoF positions symbolically just stand for the number of basis functions on each element and not for the corresponding position of nodal interpolation points.
\end{remark}

Now, we define the cell-based Vanka smoother for the proposed space-time finite element approximation. This part generalizes previous work on the local Vanka smoother to the higher order time discretization. On the mesh level $g$ (cf.\ Fig.~\ref{Fig:MG_overview}) and in a single iteration step, the cell-based Vanka smoother is applied to all $I_n$-coefficient subvectors of the spatial degrees of freedom corresponding to the respective mesh cell $T$. On grid level $g$, with $g\in \{2,\ldots,G\}$, we let the solution vector $\vec d$ of \eqref{eq:LinSysMG} be subdivided in terms of velocity and pressure subvectors according to the structure of the solution vector $\vec X_n$, defined in \eqref{eq:totdofs}, of the nonlinear system \eqref{eq:nonlineq}, such that 
\begin{equation}
	\label{eq:totdofs_2}
	\vec d = (\vec {d}^{0}_{1},\ldots,\vec {d}^{0}_{d},\vec {q}^{0},\ldots, \vec {d}^{k}_{1},\ldots,\vec {d}^{k}_{d},\vec {q}^{k})^\top \in \R^{(k+1)\times (d\cdot R_g + S_g)}\,.
\end{equation}
Here, $R_g$ and $S_g$ denote the number of (global) degrees of freedom for the velocity and pressure variable on grid level $g$, where $\vec X_n$ is defined for the finest mesh level $\mathcal T_G$. The subvectors $\vec {d}^{l}_{1},\ldots , \vec {d}^{l}_{d}$ for $l=0,\ldots, k$, correspond to the velocity values (or their corrections, respectively) and the subvectors $\vec {q}^{l}$ for $l=0,\ldots, k$, to the pressure values (or their corrections, respectively) on the grid level $g$.
With the amount $n_p$ of the local pressure degrees of freedom on each element, $n_p = \binom{d+r}{r}$, we then denote by $\vec d_T$ the subvector of $\vec d$ that is built from the degrees of freedom in $\vec d$ that are associated with the element $T$, such that 
\begin{equation}
	\label{eq:locdofs}
	\vec d_T = (\vec {d}^{0}_{1,T},\ldots,\vec {d}^{0}_{d,T},\vec {q}^{0}_T,\ldots, \vec {d}^{k}_{1,T},\ldots,\vec {d}^{k}_{d,T},\vec {q}^{k}_T)^\top \in \R^{(k+1)\times (d\cdot (r+1)^d + n_p)}\,.
\end{equation}
Here, the subvectors $\vec {d}^{l}_{1,T},\ldots, \vec {d}^{l}_{d,T}$ for $l=0,\ldots, k$, correspond to the velocity values on the element $T$ and the subvectors $\vec {q}^{l}_T$ for $l=0,\ldots, k$, to the pressure values. Further, for right-hand side vector $\vec r$ of \eqref{eq:LinSysMG} on grid level $g$ we let 
\begin{equation}
	\label{eq:defr0}
	\vec r_0 = (\vec {r}^{0}_{1},\ldots,\vec {r}^{0}_{d},\vec 0^0,\ldots, \vec {r}^{k}_{1},\ldots,\vec {r}^{k}_{d},\vec {0}^{k})^\top \in \R^{(k+1)\times (d\cdot R_g + S_g)}
\end{equation}
with the partition of $\vec r$ into subvectors induced by \eqref{eq:totdofs_2}.
On a single cell $T$, the local Jacobian matrix $\mat{J}_T$ is defined as
\begin{equation*}
	\vec J_{T}  = \dfrac{\partial \vec F_{T}}{\partial \vec X_T}(\vec X_T) \,,
\end{equation*}
where, in contrast to \eqref{eq:definition_jacobian_global}, we skipped the index of the time interval $n$ and the Newton step $m$ for brevity.
On such a cell $T$, the smoothing operator $\vec S_T(\vec{d}, \vec{r})$ is then defined by 
\begin{equation}
\label{eq:LocVanka}
\vec S_T(\vec{d}, \vec{r}) := \vec d_T + \vec J_T^{-1} (\vec r_0 - \vec J \vec d)_T
\end{equation}
In \eqref{eq:LocVanka}, the vector $(\vec r_0 - \vec J \vec d)_T\in \R^{(k+1)\times (d\cdot (r+1)^d + n_p)}$ denotes the local subvector of $(\vec r_0 - \vec J \vec d)\in \R^{(k+1)\times (d\cdot R_l + S_l)}$, that is obtained by condensing $(\vec r_0 - \vec J \vec d)$ to its components corresponding to the mesh cell $T$, similarly to \eqref{eq:locdofs}.
\textblue{The global vector is here computed fully in parallel, utilizing Trilinos functions.}
The full application of the smoother $\vec S(\vec{d}, \vec{r})$ then comes through running over all cells of the corresponding mesh level and applying the local smoother $\vec S_T(\vec{d}, \vec{r})$ to each of the elements \textblue{by an updating strategy, similar to the Jacobi iteration method.
Since the finite elements for the velocity are continuous, each velocity degree of freedom on interior faces of a cell is updated at least twice.
In our implementation we simply overwrite the corresponding values of the vector $\vec{d}$, that is being smoothed.
Thus, a degree of freedom connected with multiple cells is determined by its last update in the loop over all cells.}

The appreciable advantage of the Vanka smoother is that the system, that has to be solved on each cell, or the inverse of the local Jacobian matrix $\mat{J}_T^{-1}$, respectively, is small compared to the global system with system matrix $\mat{J}$.
This will be addressed further in the next subsection. The efficiency of the application of the Vanka smoother in complex simulations with a high number of mesh cells depends on two ingredients:
\begin{enumerate}
	\item The efficient application of $\mat{J}_T^{-1}$: How are the local systems defined by \eqref{eq:LocVanka} solved? 
	
	\item The efficient data exchange in the parallel environment: How are the data for computing $\mat{J}_T$ or $\mat{J}_T^{-1}$, respectively, assembled? 
\end{enumerate}
These two issues are discussed in the following. 

\subsection{Efficient application of $\mat{J}_T^{-1}$}

The implementation of the operator $\mat{J}_T^{-1}$ is an important ingredient for the efficiency of the GMG approach in computations, since the Vanka smoother is applied . We recall that the GMG method is used as a preconditioner in GMRES iterations for solving the Newton linearized system of each subinterval $I_n$. We also refer to Fig.~\ref{Fig:MG_overview} illustrating the usage of smoothing steps on the grid levels of a GMG V-cycle.
In our implementation of the GMG method, inverses of the element-wise Jacobian matrices $\mat{J}_T^{-1}$, for all $T\in \mathcal T_{h_l}$ with $l=1,\ldots,L$, are pre-computed after each update of the Jacobian matrix $\mat J$. For this, we use LAPACK routines to pre-compute the matrices $\mat{J}_T^{-1}$ and store them in a hashed \emph{unordered\_map}. If $\mat{J}_T^{-1}$ has to be applied on a cell $T$ according to \eqref{eq:LocVanka}, the costs for looking up the corresponding inverse is an operation with a complexity of order $\mathcal{O}(1)$.
\\
The costs in terms of memory for storing each inverse $\mat{J}_T^{-1}$ is, for instance, $85 \cdot 85 \cdot \SI{64}{\bit} = \SI{57.8}{\kilo\byte}$ or $\SI{0.0578}{\mega\byte}$ for a three-dimensional problem for the spatial approximation by the $\mathbb Q_2$--$\mathbb P_{1}^{\text{disc}}$ pair of finite elements (corresponding to $r=2$ in \eqref{Def:Vh0}) and the $dG(0)$ time discretization (corresponding to $k=0$ in \Cref{eq:time_discrete_ansatz}) on a 64-bit machine, plus some (negligible) additional overhead to store, for instance, the hashes. For a $dG(1)$ time discretization (corresponding to $k=1$ in \Cref{eq:time_discrete_ansatz}), the local Jacobian $\mat{J}_T$ is a $170 \times 170$ matrix and the needed amount of memory is $170 \cdot 170 \cdot \SI{64}{\bit} = \SI{231.2}{\kilo\byte}$ or $\SI{0.231}{\mega\byte}$.
Since the code is parallelized, every process has to store only the information, data and inverses of the cells that it owns. Therefore, the additional amount of memory, that is needed in each process, can be decreased by increasing the number of involved processors.

\subsection{Efficient data exchange in parallel environments}

\begin{figure}[h!tb]
	\centering
	\includegraphics[width=0.6\textwidth,keepaspectratio]
	{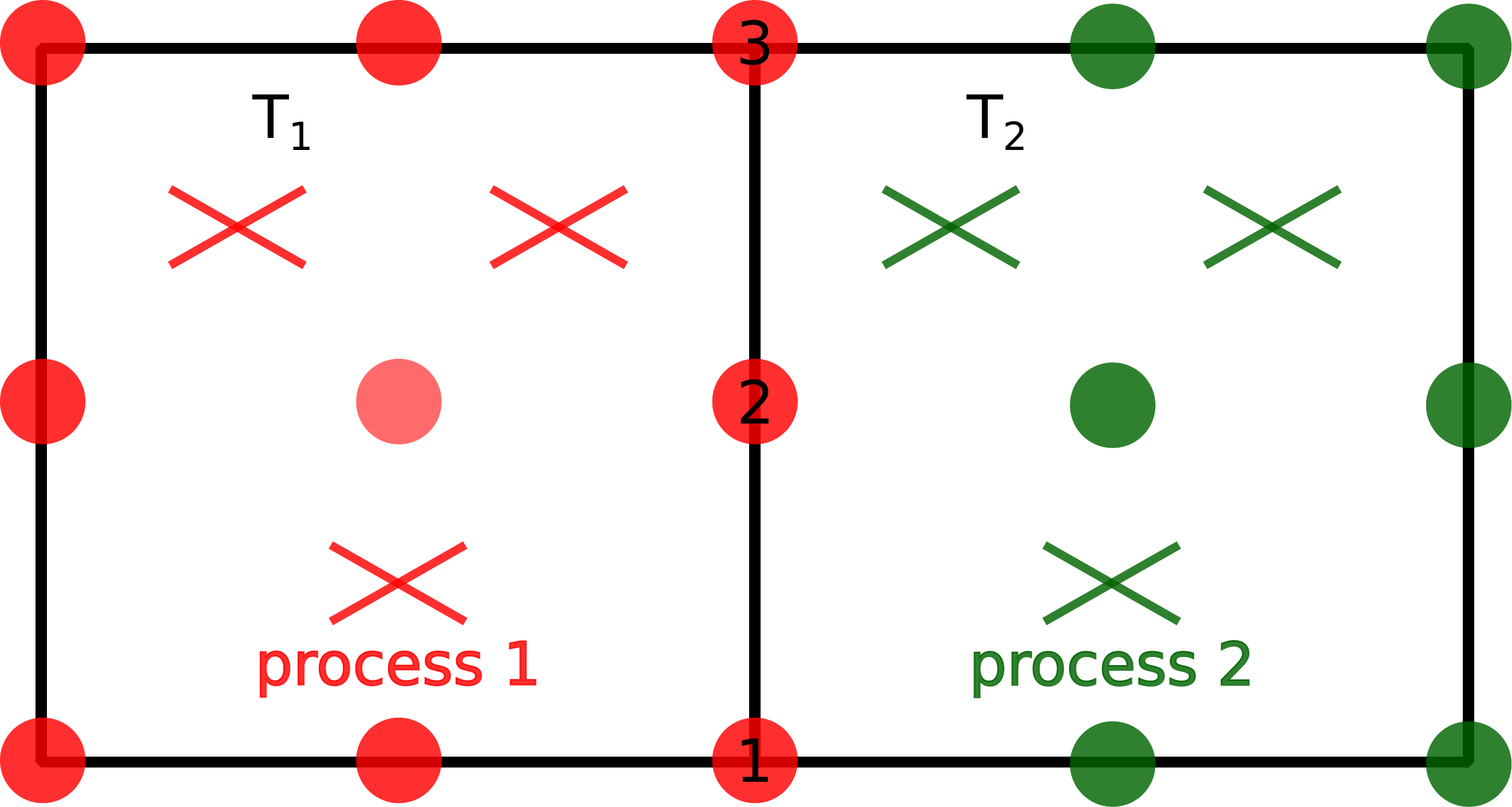}
	\caption{Two cells that belong to two different processes with the DoFs for $\mathbb{Q}_2$--$\mathbb{P}_{1}^{\text{disc}}$ elements}
	\label{Fig:Parallel_Problem}
\end{figure}

For pre-computing the inverses of the local Jacobian matrices  $\mat{J}_T^{-1}$, the entries $\mat{J}_{i,j}$ of $\mat{J}_T$ in the element $T$ have to be computed. If the code is executed in parallel by multiple processes, the data access problem that is sketched in Fig.~\ref{Fig:Parallel_Problem} occurs. When the local matrix $\mat{J}_{T_1}$ on $T_1$ is assembled by the process 1, all the needed matrix entries of the global  Jacobian matrix $\mat{J}$ are available, and can be copied to the local Jacobian $\mat{J}_{T_1}$, since process 1 owns all the involved degrees of freedom. In contrast, process 2 doesn't own the degrees of freedom on the face separating $T_1$ and $T_2$, since in a parallel environment every process has only read access to the entries it owns. For computing $\mat{J}_{T_2}$, the entries in the global matrix $\mat{J}$ of process 1, corresponding to the degrees of freedom on the common interface, are required.

To provide and exchange the needed data efficiently, the following data structure, called \emph{map\_proc\_row\_column\_value}, is generated on each of the involved processes. It involves the following, nested containers (from top to bottom).
\begin{itemize}
	\item \emph{dealii::MGLevelObject\qdist{\,\,}}: The top level object, that contains the next elements for each mesh level $l$. The \emph{MGLevelObject\qdist{\,\,}} is basically a container like \emph{std::vector\qdist{\,\,}}, but with the option to shift indexing. So if the coarsest mesh starts e.\,g. on level $l = 2$ one can access the elements inside this object with the $[\,]$ operator and an index, starting from 2 onward.
	\item \emph{std::map\qdist{unsigned int, std::unordered\_map\qdist{$\ldots$} }}: A map, whose keys are the process numbers (an unsigned int) of the neighboring processes, that own certain degrees of freedom. These are exactly those degrees of freedom, that the owning process needs to access during the assembly of the local Jacobian $\mat{J}_T$.
	The process numbers are obtained by \Cref{alg:multigrid_parallel}.
	The value of the map is a (hashed) unordered\_map, which leads to the next container inside the structure:
	\item \emph{std::unordered\_map\qdist{ std::pair\qdist{unsigned int, unsigned int}, double}}: For each neighboring process a hashed, unordered\_map is stored, that contains the global row and column number (both unsigned int) of the needed matrix entries and assigns them to the corresponding value, which is stored as a double.
	\textblue{Internally the std::pair of global row and column numbers is stored as hashed value and therefore accessing or inserting operations into this data structure have an average complexity of $\mathcal{O}(1)$.}
\end{itemize}
After generation of the mesh hierarchy, every process executes the Algorithm \ref{alg:multigrid_parallel}.


\begin{algorithm}[H]
	\SetAlgoLined
	Create \emph{map\_proc\_row\_column\_value} \;
	\ForEach{mesh level $l$}{
		\ForEach{locally owned cell $K_i$}{
			\ForEach{DoF $j$}{
				\If{DoF $j$ is not owned by this process}{
					Get the number of the process $np$ that owns it \;
					Add $np$ and all couplings of DoF $j$ with all other locally owned DoFs to \emph{map\_proc\_row\_column\_value} and assign it a value of 0. \;	
	}}}}
	\caption{Determine values of $\mat{J}$ owned by foreign processes}
	\label{alg:multigrid_parallel}
\end{algorithm}

\begin{remark}
	If the underlying mesh or the distribution of the degrees of freedom to the involved processes is fixed, which is especially the case if no remeshing is necessary between time-steps, then \Cref{alg:multigrid_parallel} needs to be executed only once in the simulation. For instance, this holds in the numerical examples of Sec.~\ref{sec:numerical_example}.
	In the simulation of flow problems on evolving domains by CutFEM approaches, that are currently focused strongly (cf.~\cite{vonwahlUnfittedEulerianFinite2020,burmanEulerianTimesteppingSchemes2020}), this applies similarly and results in negligible computational costs.
	\\
\textblue{The number of the owning process of a not locally owned DoF is computed using the function \emph{compute\_index\_owner()} of the \emph{dealii::Utilities::MPI} namespace, which uses non-blocking point-to-point communication.}
	
\end{remark}
Afterwards the simulation is continued until all contributions of the global system matrix $\mat{J}$ on all levels $g$ are assembled.
For building the Vanka smoother by assembling and storing the local matrices $\mat{J}_T^{-1}$, the respective sparse matrix entries have to be exchanged such that every process can access the entries of each local Jacobian matrix $\mat{J}_T$. This is done by Algorithm~\ref{alg:multigrid_parallel_1}.

\begin{algorithm}[H]
	\SetAlgoLined
	\ForEach{mesh level $g$}{
		Use MPI some-to-some communication to transfer information: \emph{recv\_proc\_row\_column\_value} $\leftarrow$ \emph{map\_proc\_row\_column\_value} \;
		\ForEach{proc $p$ in \emph{recv\_map\_proc\_row\_column\_value}}{
			\ForEach{global row $i$ and column $j$ in \emph{recv\_proc\_row\_column\_value}}{
				value $\leftarrow$ $\mat{J}_{i,j}$
		}}
		Use MPI some-to-some communication to transfer back information: \emph{map\_proc\_row\_column\_value} $\leftarrow$ \emph{recv\_proc\_row\_column\_value} \;
	}
	\caption{Update \emph{map\_proc\_row\_column\_value} on each process}
	\label{alg:multigrid_parallel_1}
\end{algorithm}

\begin{remark}
	\Cref{alg:multigrid_parallel_1} has to be executed whenever the global Jacobian $\mat{J}$ in \eqref{eq:LinSysMG} is updated.
	To block as less resources as necessary, just the processes that need to exchange information, communicate with each other.
	The exchange of data is implemented in the object  \emph{map\_proc\_row\_column\_value}.
	The data transfer is done entirely in one single step (see \Cref{alg:multigrid_parallel_1}), instead of querying the data, that is needed for a single cell $T$ from foreign processes, in the assembly routine.
	\\
	\textblue{Each process saves the received information in the temporary object \emph{recv\_map\_proc\_row\_column\_value}, which is a container of the same type as \emph{map\_proc\_row\_column\_value}, but the process numbers are the one of the processes, querying the information.
	After receiving this object, each process looks up the queried matrix values in the corresponding rows and columns.
	In the last step this information is transferred back to the querying processes and stored in \emph{recv\_map\_proc\_row\_column\_value}.}
	By this approach, we reduce the communication between processes to an absolute necessary minimum.
	\\
	\textblue{The some-to-some communication utilizes the \emph{some\_to\_some()} function of the \emph{dealii::Utilities::MPI} namespace, which basically relies on two-sided communication, utilizing \emph{MPI\_Isend()} to exchange the information.
	Before the actual communication, the number of processes, that send information tho this process, is computed using a \emph{MPI\_Reduce\_scatter\_block()} call, so each process knows in advance how many calls to \emph{MPI\_Recv()} are necessary.
	The actual sender is determined using \emph{MPI\_Probe()} functionality.
}
\end{remark}

\begin{remark}
	To call MPI functions that transfer data between processes, these data have to be serializable. The C++17 version of \emph{std::unordered\_map} included in the standard library is by default not serializable.
	In the code of this work, the boost C++ library is used that provides serialization capabilities for \emph{std::unordered\_map} via the interfaces \emph{serialization.hpp} and \emph{unordered\_map.hpp}.
\end{remark}

\section{Numerical examples}
\label{sec:numerical_example}

In the following we analyze computationally the performance properties of the proposed GMG approach. This is down for the well-known benchmark problems of flow around a cylinder; cf.~\cite{schaferBenchmarkComputationsLaminar1996}. Our computations were done on a Linux cluster with 96 nodes, each of them with 2 CPUs and 14 cores per CPU. The CPUs are \emph{Intel Xeon E5-2680 v4} with a base frequency of $\SI{2.4}{\giga\hertz}$, a maximum turbo frequency of $\SI{3.3}{\giga\hertz}$ and a level 3 cache of $\SI{35}{\mega\byte}$. Each node has $\SI{252}{\giga\byte}$ of main memory. In this work, scaling experiments on up to the user limit of 32 nodes were performed.

\subsection{Flow around a cylinder in two space dimensions}
\label{sec:dfg_2d}

\begin{figure}[h!tb]
	\centering
	\includegraphics[width=0.9\linewidth]{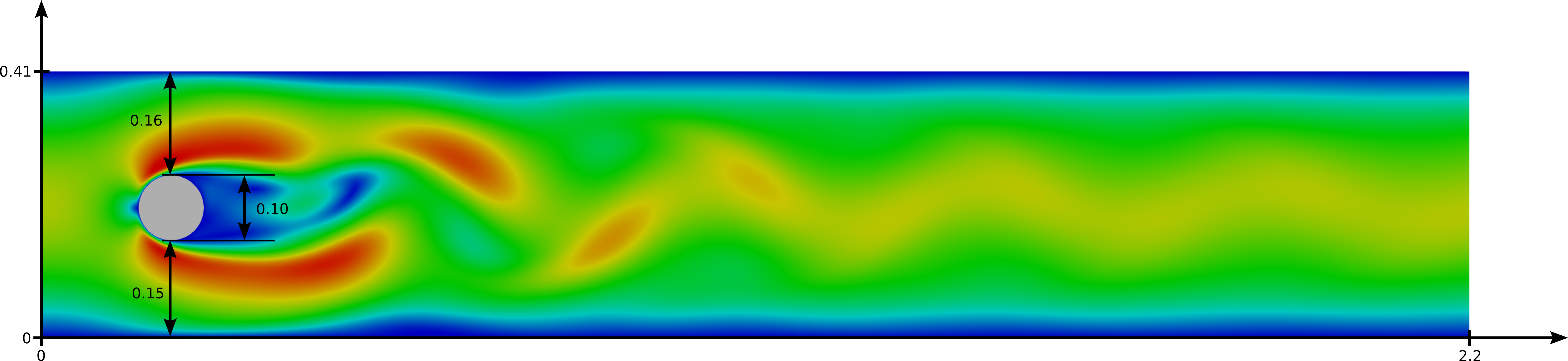}
	\caption{Geometrical setting of the 2d DFG benchmark with computed velocity profile of the fully developed flow.}%
	\label{fig:problem_overview_2D_DFG}
\end{figure}
In the first numerical experiment we consider the well-known 2d DFG benchmark setting of flow around a cylinder, defined in \cite{schaferBenchmarkComputationsLaminar1996}. The problem setting is illustrated in Fig.~\ref{fig:problem_overview_2D_DFG}. Quantities of interest and comparison in the simulations are the drag and lift coefficient of the flow on the circular cross-section (cf.~\cite{schaferBenchmarkComputationsLaminar1996}). With the drag and lift forces $F_D$ and $F_L$ on the rigid circle $S$ given by 
\begin{align}
	\label{eq:drag_lift_forces}
	F_D &= \int_S
	\left(
	\nu \frac{\partial v_t}{\partial \vec{n}} n_y - P n_x
	\right) \d S \,,
	&
	F_L &= - \int_S
	\left(
	\nu \frac{\partial v_t}{\partial \vec{n}} n_x - P n_y
	\right) \d S \,,
\end{align}
where $\vec n$ is the normal vector on $S$ and $v_t$ is the tangential velocity $\vec{t} = (n_y, -n_x)^\top$, the drag and lift coefficient $c_D, c_L$ are defined by means of
\begin{align}
	c_D &= \frac{2}{\bar{U}^2 L} F_D \,,
	&
	c_L &= \frac{2}{\bar{U}^2 L} F_L \,.
\end{align}
According to \cite{schaferBenchmarkComputationsLaminar1996}, we choose the viscosity $\nu = 0.001$ and the boundary condition on the inflow boundary $\Gamma_i$ as $\vec{g}_i(x, y, t) = ({4 \cdot 1.5 \cdot y (0.41 - y)}/{0.41^2}, 0)^\top$. This leads to a Reynolds number of $Re = 100$ and a time-periodic flow behavior. The final simulation time is $T=10$, such that $I = (0, 10]$. The space--time discretization is done in the spaces $(\mathbb P_1(I_n;H_h^2))^2 \times \mathbb P_1(I_n;H_{h,\text{disc}}^1)$ with time step size $\tau = \num{0.005}$.
The Newton iteration is stopped when the residual of the nonlinear equation is smaller than $\num{1e-10}$ or a relative reduction of the initial residual by a factor of $\num{10000}$ is reached.
The computations were done on 4 nodes of the Linux cluster. The computed velocity field of the fully developed flow is presented in Fig.~\ref{fig:problem_overview_2D_DFG}.

Table~\ref{tab:2d_DFG_results} shows the space-time degrees of freedom in one single time step and summarizes for three mesh levels the computed maximum drag and lift coefficients of the fully developed flow as well as the average number of Newton iterations per time step and the average number of GMRES iterations per Newton step.
\textblue{The coarse level is always set to $g = 1$ and we used 8 - Nr. mesh levels. So for instance for simulation Nr. 1 the finest mesh level is level 6.}
Table~\ref{tab:2d_DFG_results} further shows that on all mesh levels an average of less than two Newton iterations per subinterval is obtained. Moreover, the number of GMRES iterations with GMG preconditioning remains (almost) grid independent. This demonstrates the efficiency of the proposed approach.  

{
	\sisetup{scientific-notation = false,
		round-mode=places,
		round-precision=4,
		output-exponent-marker=\ensuremath{\mathrm{e}},
		table-figures-integer=1, 
		table-figures-decimal=3, 
		table-figures-exponent=1, 
		table-sign-mantissa = false, 
		table-sign-exponent = true, 
		table-number-alignment=center} 
	
	\begin{table}[h!t]
		\caption{Computed drag and lift coefficients and average number of Newton steps per time step and of GMRES iterations per Newton step in the two-dimensional benchmark \textblue{for three simulations.$h_{\text{max}}$ is the maximum diameter of the cells on the finest mesh level $G$ and DoFs denote the problem size on this level.} Reference values: $c_{D_{max}} \in [\num{3.2200}, \num{3.2400}]$ and $c_{L_{max}} \in [\num{0.9900}, \num{1.0100}]$, see \cite{schaferBenchmarkComputationsLaminar1996}.}
		\centering
		\begin{tabular}{c@{\hskip 4ex} c@{\hskip 4ex} b@{\hskip 4ex} c@{\hskip 4ex}  c@{\hskip 4ex} c@{\hskip 4ex} c@{\hskip 4ex} b@{\hskip 4ex}}
	\toprule
	{Nr.} & {DoFs} & {$h_{\text{max}}$} &{$c_{D_{max}}$} & {$c_{L_{max}}$} & {$\bar{n}_{\text{Newton}}$} & {$\bar{n}_{\text{GMRES}}$}  & {$\bar{n}_{\text{GMRES}}$} \\
	\cmidrule(lr){1-6} \cmidrule(r){7-7} \cmidrule(r){8-8}
	{1} & \num{8305664} & \num{0.0022} & \num{3.2350} & \num{1.0062} & 1.72 & \num{10} & \num{7} \\
	{2} & \num{2080256} & \num{0.0048} & \num{3.2227} & \num{1.0060} & 1.67 & \num{9} & \num{7} \\
	{3} & \num{521984} & \num{0.0090} &\num{3.1274} &  \num{0.9637} & 1.53 & \num{9} & \num{6}\\
	\cmidrule(lr){1-6} \cmidrule(r){7-7} \cmidrule(r){8-8}
	\multicolumn{6}{c}{Multigrid cycle} & {$V(1,1)$} & {$V(4,4)$} \\
	\bottomrule
\end{tabular}
		\label{tab:2d_DFG_results}
	\end{table}
}

Table~\ref{tab:2d_DFG_timings} summarizes the wall time consumed by the different parts of the algorithms.
$t_{\text{GMRES}}$ summarizes the time, spent in the outer flexible GMRES solver (including the GMG preconditioner).
$t_{\text{Inv}}$ is the time spent for computing the local inverses of the cell Jacobians $J_T^{-1}$ and $t_{\text{Upd}}$ is the time that is spent in \Cref{alg:multigrid_parallel_1}, exchanging necessary information between processes.
\textblue{The rest of the computation time is spent mainly in thee assembly routine for the system Jacobi matrix and to a negligible amount in data output routines.}
\\
In our simulations, the same number of compute nodes (4 nodes) were used \textblue{in all simulations and on all multigrid levels}. The usage of 4 nodes for the considered problem dimensions leads to a great difference in the percent wall time of the GMRES solver. \textblue{In problem setting Nr. 1}, the GMG preconditioned GMRES iterations consumed only $\SI{17.96}{\percent}$ of the total wall time.
The latter indicates that the benefit of a faster system assembly, by using more nodes, would have paid off and annihilated the costs of an increased parallel communication.
In contrast, on the coarsest level the number of nodes was set too high to pay off. We observe just a slight decrease in the wall time despite nearly quartering the number of degrees of freedom.

{
	\sisetup{scientific-notation = false,
		round-mode=places,
		round-precision=2,
		output-exponent-marker=\ensuremath{\mathrm{e}},
		table-figures-integer=1, 
		table-figures-decimal=3, 
		table-figures-exponent=1, 
		table-sign-mantissa = false, 
		table-sign-exponent = true, 
		table-number-alignment=center} 
	
	\begin{table}[h!t]
		\small
		\caption{Wall time consumption of the two-dimensional benchmark simulation.}
		\label{tab:2d_DFG_timings}
		\begin{subtable}[t]{\textwidth}
		\centering
		\caption{Utilizing a $V(1,1)$ multigrid cycle.}
		\label{tab:2d_DFG_timings_a}
	\begin{tabular}{c@{\,\,\,\,}c  c@{\,\,\,\,}c   c@{\,\,\,\,}c  c@{\,\,\,\,}c}
		\toprule
		{DoFs} & {$t_{\text{wall}}$} &
		{ $t_{\text{GMRES}}$ } & {\% of $t_{wall}$} &
		{$t_{\text{Inv}}$} & {\% of $t_{wall}$} &
		{$t_{\text{Upd}}$} & {\% of $t_{wall}$} \\
		\cmidrule(lr){1-2}
		\cmidrule(lr){3-4}
		\cmidrule(lr){5-6}
		\cmidrule(lr){7-8}
		\num{8305664} & \SI{6.14}{\hour} &
		\SI{1.10}{\hour} & \num{17.96} &
		\SI{0.15}{\hour} & \num{2.42} &
		\SI{0.04}{\hour} & \num{0.66} \\
		\num{2080256} & \SI{1.10}{\hour} &
		\SI{0.57}{\hour} & \num{51.81} &
		\SI{0.04}{\hour} & \num{2.16} &
		\SI{0.02}{\hour} & \num{1.11} \\
		\num{521984} & \SI{0.83}{\hour} &
		\SI{0.4}{\hour} & \num{48.16} &
		\SI{0.01}{\hour} & \num{1.51} &
		\SI{0.01}{\hour} & \num{1.40} \\
		\bottomrule
	\end{tabular}
\end{subtable}
\begin{subtable}[t]{\textwidth}
	\vspace*{3ex}
	\centering
	\caption{\textblue{Utilizing a $V(4,4)$ multigrid cycle.}}
	\label{tab:2d_DFG_timings_b}
	\textblue{
	\begin{tabular}{c@{\,\,\,\,}c  c@{\,\,\,\,}c   c@{\,\,\,\,}c  c@{\,\,\,\,}c}
		\toprule
		{DoFs} & {$t_{\text{wall}}$} &
		{ $t_{\text{GMRES}}$ } & {\% of $t_{wall}$} &
		{$t_{\text{Inv}}$} & {\% of $t_{wall}$} &
		{$t_{\text{Upd}}$} & {\% of $t_{wall}$} \\
		\cmidrule(lr){1-2}
		\cmidrule(lr){3-4}
		\cmidrule(lr){5-6}
		\cmidrule(lr){7-8}
		\num{8305664} & \SI{5.76}{\hour} &
		\SI{0.72}{\hour} & \num{12.44} &
		\SI{0.15}{\hour} & \num{2.61} &
		\SI{0.04}{\hour} & \num{0.69} \\
		\num{2080256} & \SI{0.88}{\hour} &
		\SI{0.35}{\hour} & \num{40.07} &
		\SI{0.04}{\hour} & \num{4.52} &
		\SI{0.02}{\hour} & \num{2.26} \\
		\num{521984} & \SI{0.70}{\hour} &
		\SI{0.27}{\hour} & \num{38.53} &
		\SI{0.01}{\hour} & \num{1.43} &
		\SI{0.01}{\hour} & \num{1.42} \\
		\bottomrule
	\end{tabular}
	}
	\end{subtable}
	\end{table}
}

\textblue{
\subsection{Parallel scaling}
\label{sec:GMG_parallel_scaling}
In this section we analyze the parallel performance properties of our code by a strong scaling benchmark.
We first define the parallel speedup $S$ of a program, according to \cite{amdahlValiditySingleProcessor1967}:
\begin{equation}
\label{eq:speedup}
	S = \frac{1}{r_s + \frac{r_p}{np}}\,,
\end{equation}
where $r_s$ is the ratio of the sequential fraction of the program and $r_p$ the portion, that can be scheduled in parallel with $np$ number of processes.
This is called Amdahl's Law.
So, if the problem size is fixed, the parallel speedup is limited by the serial part of the code.
Amdahl's law is under the assumption of an instant communication over a network, infinitely fast.
In practice this is not possible and therefore one also has to consider the communication costs, that are introduced when increasing the number of processes $np$, due to the finite bandwidth and latency of the network.
In practice, the speedup $S$ for a simulation, that is run multiple times on a different amount of nodes $n$, is calculated with
\begin{equation*}
	S = \frac{t_{\text{wall}}(n = n_{\text{min}})}{t_{\text{wall}}(n)}\,.
\end{equation*}
Here $n_{\text{min}}$ is the simulation with the smallest amount of nodes.

To measure the speedup S of our code, we use the spatial setup of the 2d DFG benchmark of \Cref{sec:dfg_2d}, but set the inflow condition on $\Gamma_i$ as
\begin{equation*}
	\vec{g}_i(x, y, t) = \left( \frac{4 \cdot 0.3 \cdot y (0.41 - y)}{0.41^2}, 0 \right)^\top \,.
\end{equation*}
With $\nu = 0.001$ This results in a Reynolds number of $Re = 20$.
After about a simulation time about $t = 2.3$ the flow is fully developed, which results in a static flow profile.
The final simulation time is put to $T=3$ such that $I = (0, 3]$ and the time step size is fixed to $0.1$.
\\
For the first benchmark, the numerical approximation is done in the space--time finite element spaces $(\mathbb P_1(I_n;H_h^2)^2\times \mathbb P_1(I_n;H_{h,\text{disc}}^1)$.
The mesh consists of \num{376832} cells, which results in \num{8305664} space--time degrees of freedom in each time interval.

\begin{figure}[!h]
	\centering
	\subcaptionbox{Benchmarked and ideal wall time. \label{fig:mpi_bench_abs_time}}
	[0.48\columnwidth]{
		\includegraphics[width=0.45\textwidth,keepaspectratio]{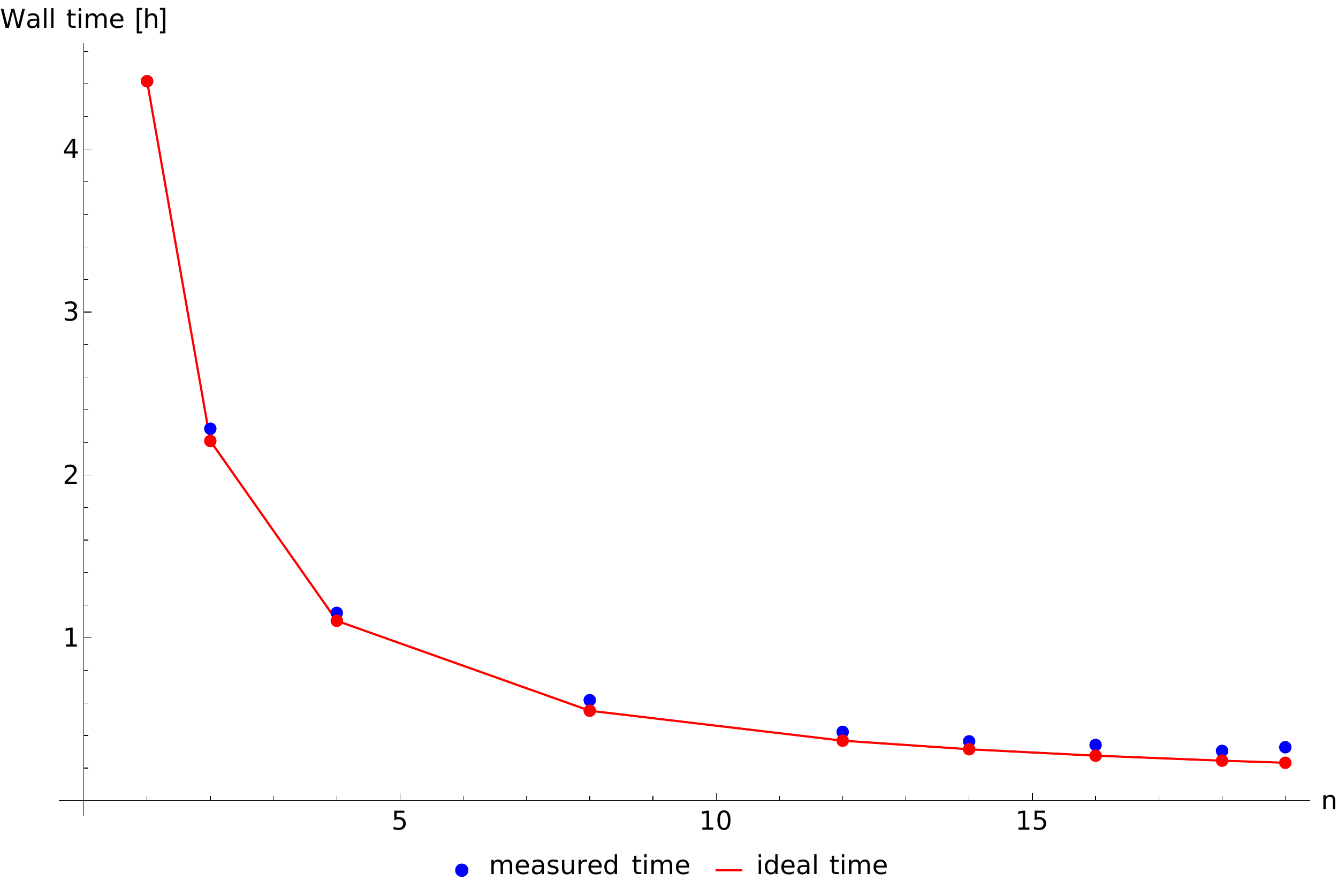}
	}
	\subcaptionbox{Benchmarked and ideal speedup $S$. \label{fig:mpi_bench_speedup}}
	[0.48\columnwidth]{
		\includegraphics[width=0.45\textwidth,keepaspectratio]{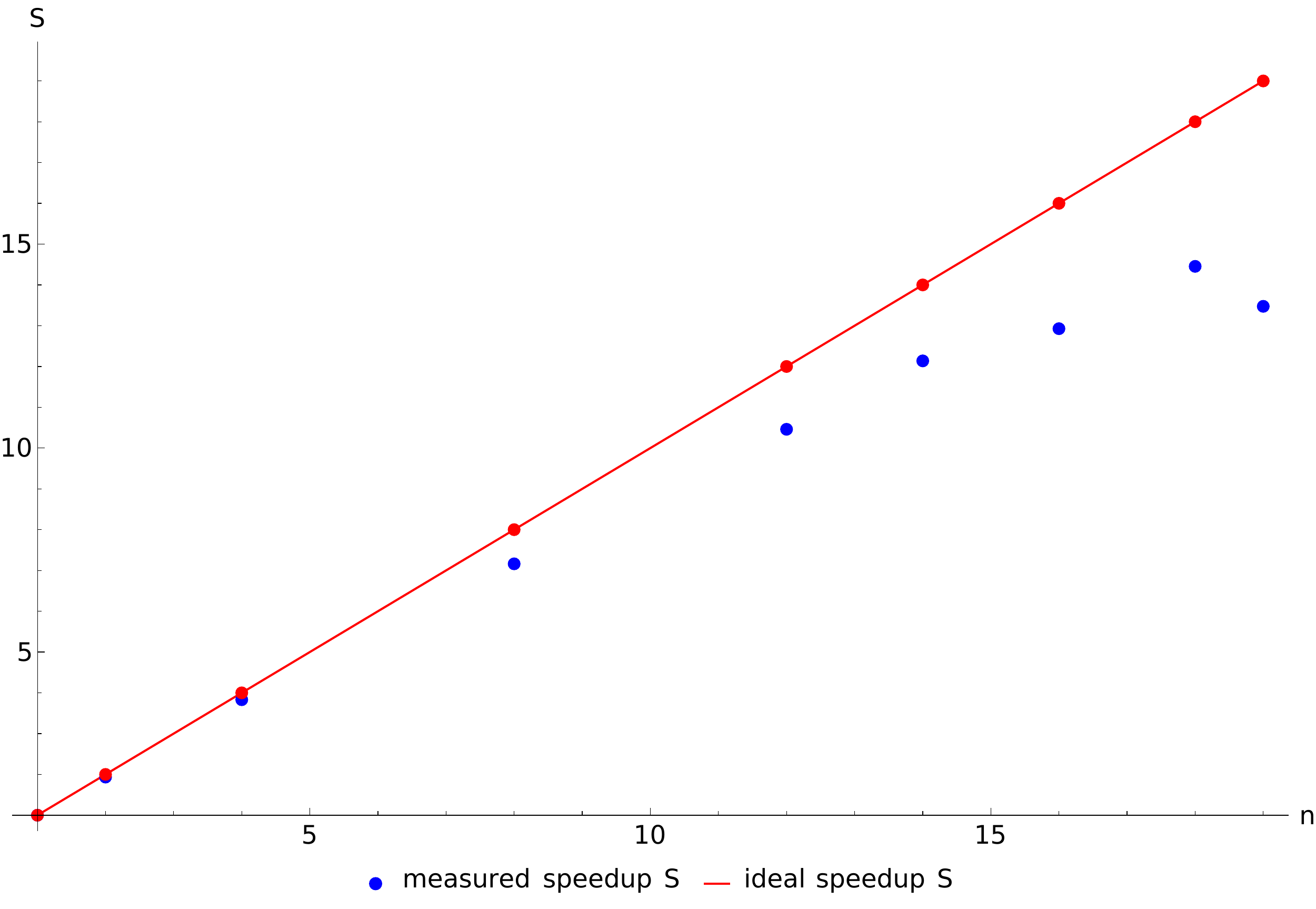}
	}
	\caption{Strong scaling results for the simulation in $(\mathbb P_1(I_n;H_h^2)^2\times \mathbb P_1(I_n;H_{h,\text{disc}}^1)$.}
	\label{fig:mpi_benchmark_results}
\end{figure}
In the benchmark the number of nodes $n$ is varied.
We assign to each node 28 processes, so that each physical CPU core owns a single process.
\Cref{fig:mpi_benchmark_results} shows the results of the benchmarks.
The ideal time is computed by setting $r_s$ in \Cref{eq:speedup} to zero.
We see nearly ideal scaling properties until the usage of 4 nodes.
Afterwards, when we further increase the number of nodes, we still see a decrease of the overall wall time of the simulation until we reach peak performance with 18 nodes.
When using even more nodes, the communication costs dominate and lead to an increase of the overall runtime of the simulation.
In this example the wall time could be reduced from \SI{4.42}{\hour} to just \SI{19.67}{\minute}, when using 18 nodes with 504 processes.

For the second scaling benchmark, the numerical approximation is done in the space--time finite element spaces $(\mathbb P_2(I_n;H_h^3)^2\times \mathbb P_2(I_n;H_{h,\text{disc}}^2)$.
We use the same spatial mesh as before, with \num{376832} cells.
That results this time in \num{27166464} space--time degrees of freedom in each time interval.
The minimum number of nodes, that were used for this benchmark, was 4.
With this configuration the runtime was \SI{54.17}{\hour}.
The maximum number of nodes, that we used, was 64.
With this configuration the runtime was reduced to \SI{3.39}{\hour}.
\Cref{fig:mpi_benchmark_results_2} shows the results of the benchmarks.
\begin{figure}[!ht]
	\centering
	\subcaptionbox{Benchmarked and ideal wall time. \label{fig:mpi_bench_abs_time_2}}
	[0.48\columnwidth]{
		\includegraphics[width=0.45\textwidth,keepaspectratio]{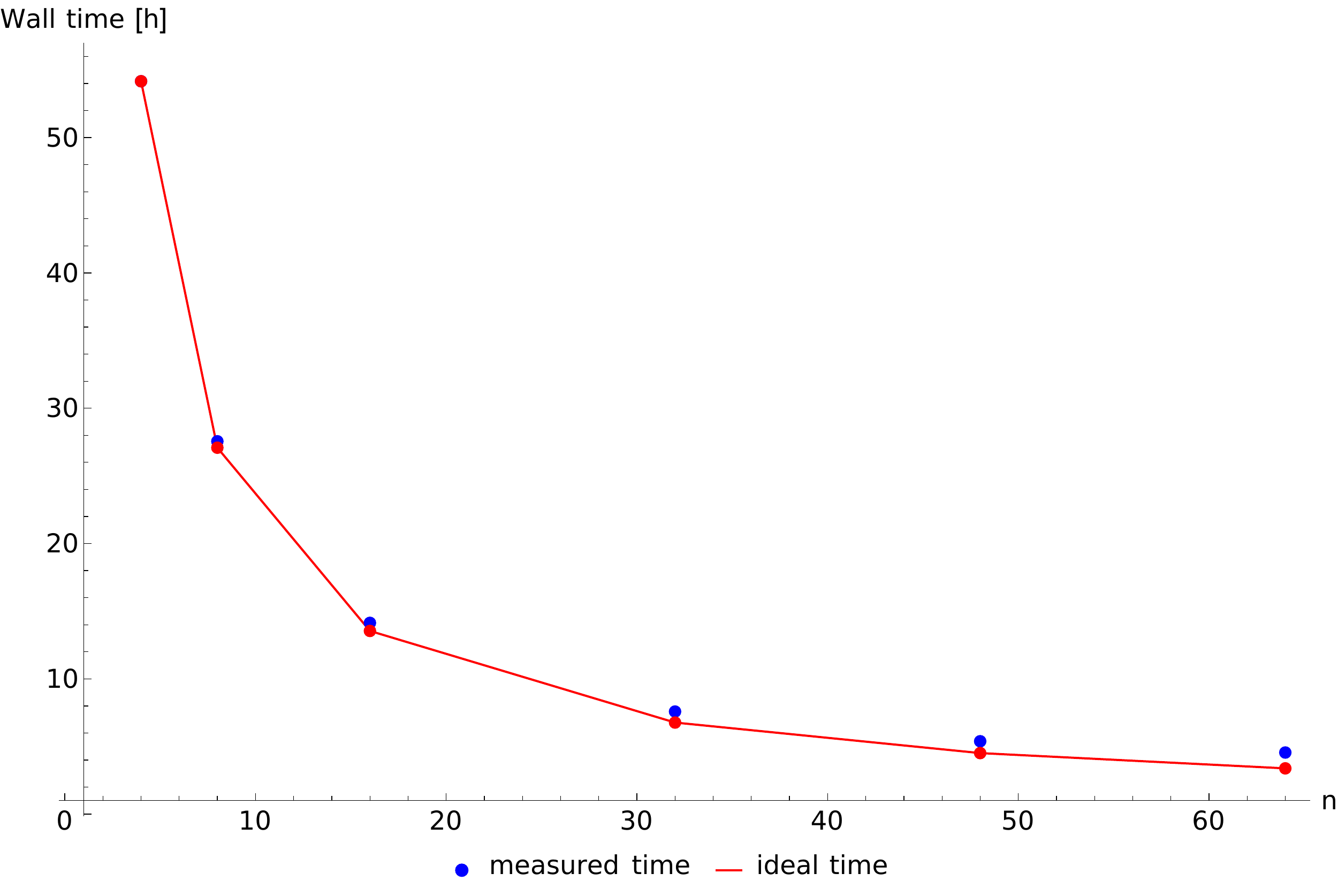}
	}
	\subcaptionbox{Benchmarked and ideal speedup $S$. \label{fig:mpi_bench_speedup_2}}
	[0.48\columnwidth]{
		\includegraphics[width=0.45\textwidth,keepaspectratio]{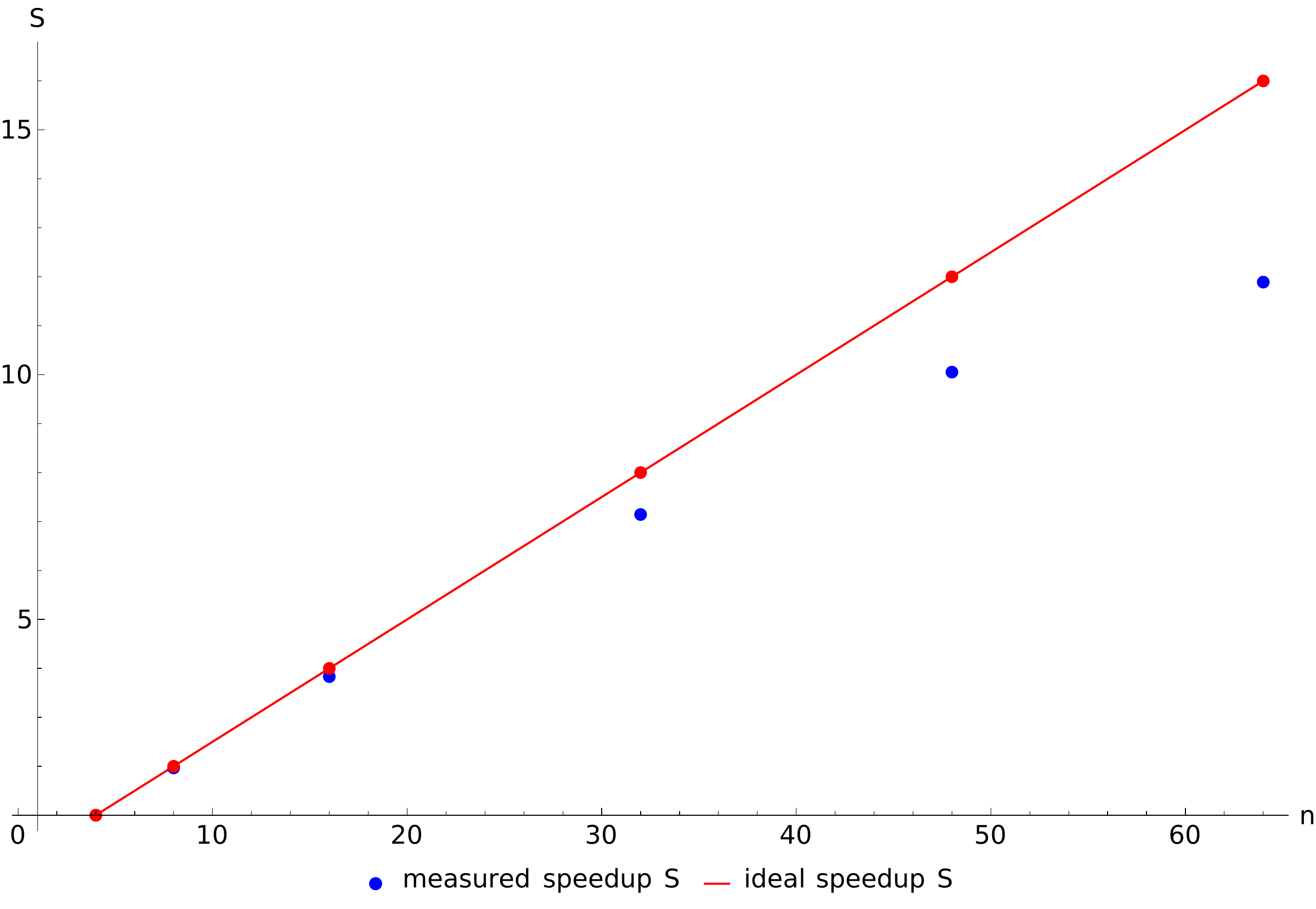}
	}
	\caption{Strong scaling results for the simulation in $(\mathbb P_2(I_n;H_h^3)^2\times \mathbb P_2(I_n;H_{h,\text{disc}}^2)$.}
	\label{fig:mpi_benchmark_results_2}
\end{figure}
Up to 16 nodes we have nearly optimal scaling results.
Then one observes an increasing difference from the ideal speedup to the measured speedup, due to the increase of the communication costs.
In \Cref{tab:strong_scaling} we summarize the characteristic statistics for the two simulations.
{
	\sisetup{scientific-notation = false,
		round-mode=places,
		round-precision=2,
		output-exponent-marker=\ensuremath{\mathrm{e}},
		table-figures-integer=1, 
		table-figures-decimal=3, 
		table-figures-exponent=1, 
		table-sign-mantissa = false, 
		table-sign-exponent = true, 
		table-number-alignment=center} 

\begin{table}[!ht]
	\small
	\centering
	\caption{Number of nodes $n$, number of corresponding processes $np$, average cells, that each process owns ($\overline{\text{cells}_\text{p}}$), wall time and speedup $S$ for the two strong scaling benchmarks.}
	\label{tab:strong_scaling}
	\begin{subtable}[t]{\textwidth}
		\centering
		\caption{Using $\mathbb P_1(I_n;H_h^2)^2\times \mathbb P_1(I_n;H_{h,\text{disc}}^1)$ elements on each time interval.}
		\label{tab:strong_scaling_a}
		\begin{tabular}{c@{\hskip 4ex} c@{\hskip 4ex} c@{\hskip 4ex}  c@{\hskip 4ex} c@{\hskip 4ex}}
			\toprule
			{$n$} & {$np$} & {$\overline{\text{cells}_\text{p}}$} & {Wall time [h]} & {Speedup $S$} \\
			\midrule
			{1} & \num{28} & \num{13460} & \num{4.42} & \num{1.00} \\
			{2} & \num{56} & \num{6728} & \num{2.28} & \num{1.93} \\
			{4} & \num{112} & \num{3364} & \num{1.15} & \num{3.83} \\
			{8} & \num{224} & \num{1684} & \num{0.62} & \num{7.16} \\
			{12} & \num{336} & \num{1120} & \num{0.42} & \num{10.46} \\
			{14} & \num{392} & \num{960} & \num{0.36} & \num{12.14} \\
			{16} & \num{448} & \num{840} & \num{0.34} & \num{16.93} \\
			{18} & \num{504} & \num{748} & \num{0.31} & \num{14.45} \\
			{19} & \num{532} & \num{708} & \num{0.32} & \num{13.47} \\
			\bottomrule
		\end{tabular}
	\end{subtable}
	\begin{subtable}[t]{\textwidth}
		\vspace*{3ex}
		\centering
		\caption{Using $\mathbb P_2(I_n;H_h^3)^2\times \mathbb P_2(I_n;H_{h,\text{disc}}^2)$ elements on each time interval.}
		\label{tab:strong_scaling_b}
		\begin{tabular}{c@{\hskip 4ex} c@{\hskip 4ex} c@{\hskip 4ex}  c@{\hskip 4ex} c@{\hskip 4ex}}
			\toprule
			{$n$} & {$np$} & {$\overline{\text{cells}_\text{p}}$} & {Wall time [h]} & {Speedup $S$} \\
			\midrule
			{4} & \num{112} & \num{3364} & \num{54.17} & \num{1.00} \\
			{8} & \num{224} & \num{1684} & \num{27.08} & \num{1.97} \\
			{16} & \num{448} & \num{840} & \num{13.54} & \num{3.83} \\
			{32} & \num{896} & \num{420} & \num{6.77} & \num{7.14} \\
			{48} & \num{1344} & \num{280} & \num{4.51} & \num{10.05} \\
			{64} & \num{1792} & \num{212} & \num{3.39} & \num{11.89} \\
			\bottomrule
		\end{tabular}
	\end{subtable}
\end{table}
}

Both strong scaling benchmarks show, that the parallelization pays off.
Especially higher order space--time elements are applicable in a reasonable amount of time only with parallelization.

\subsection{Parameter robustness regarding $\nu$}
In this subsection we computationally analyze the robustness of the GMG preconditioned GMRES solver regarding changes in the fluid viscosity $\nu$.
We consider again the 2d DFG benchmark of \Cref{sec:numerical_example} in the setting Nr. 1 of \Cref{tab:2d_DFG_results}, with \num{8305664} DoFs per time interval but vary $\nu$ in our simulations.
We fix the time step size to $\tau = 0.005$ and utilize 4 pre- and post-smoothing steps on each multigrid level.
In our computational experiments we made the experience that adding a numerical damping factor $\omega_d$ to the smoother can increase the robustness of the GMG scheme:
\begin{equation*}
	\vec S_T(\vec{d}, \vec{r}) := \vec d_T + \omega_d \cdot \vec J_T^{-1} (\vec r_0 - \vec J \vec d)_T
\end{equation*}
Setting $\omega_d$ to a value of 0.7 leads to a remarkable reduction of the iteration numbers in our simulations.
\Cref{tab:2d_DFG_nu} shows the results for different simulations.
{
	\sisetup{scientific-notation = false,
		round-mode=places,
		round-precision=4,
		output-exponent-marker=\ensuremath{\mathrm{e}},
		table-figures-integer=1, 
		table-figures-decimal=4, 
		table-figures-exponent=0, 
		table-sign-mantissa = false, 
		table-sign-exponent = false, 
		table-number-alignment=center} 

\begin{table}[!ht]
	\caption{Average number of Newton and GMRES steps for varying $\nu$ and two different damping parameters $\omega_d$. $\infty$ means that the GMRES solver didn't converge within \num{1000} steps.}
	\centering
	\begin{tabular}{S@{\hskip 2ex} c@{\hskip 4ex} c@{\hskip 4ex} c@{\hskip 4ex} c@{\hskip 4ex}}
		\toprule
		{$\nu$} & {$\bar{n}_{\text{Newton}}$} & {$\bar{n}_{\text{GMRES}}$} & {$\bar{n}_{\text{Newton}}$} & {$\bar{n}_{\text{GMRES}}$} \\
		\cmidrule(r){1-1} \cmidrule(r){2-3} \cmidrule(r){4-5}
		0.001 & 1.72 & \num{10} & 1.69 & \num{3} \\
		0.0005 & 1.74 & \num{17}  & 1.72 & \num{4} \\
		0.0003 & 1.70 & \num{56}  & 1.70 & \num{11} \\
		0.0002 & {--} & $\infty$ & 1.75 & \num{47} \\
		0.0001 & {--} & $\infty$ & 1.72 & \num{65} \\
		\cmidrule(r){1-1} \cmidrule(r){2-3} \cmidrule(r){4-5}
		{$\omega_d$} & \multicolumn{2}{c}{1} & \multicolumn{2}{c}{0.7} \\
		\bottomrule
	\end{tabular}
	\label{tab:2d_DFG_nu}
\end{table}
}
Without damping, the simulation aborted at $\nu = 0.0002$, which is equivalent to a Reynolds number of 500, due to high iteration numbers (> \num{1000}) in the linear GMRES solver.
With $\omega_d = 0.7$, the GMG solver showed reasonable performance up to $\nu = 0.0001$, which is equivalent to a Reynolds number of 1000.
\begin{remark}
	We note that we didn't apply any fluid stabilization, therefore the occurrence of instabilities in our numerical scheme are expected in convection-dominated settings.
	There exist various stabilization techniques like the Streamline Upwind Petrov Galerkin (SUPG) or Flux-Correction methods \cite{richterFluidstructureInteractionsModels2017} to overcome this issue, which are out of the scope of this work.
\end{remark}

}

\subsection{Flow around a cylinder in three space dimensions}
\label{sec:dfg_3d}

\begin{figure}[h!tb]
	\centering
	\includegraphics[width=0.6\textwidth,keepaspectratio]
	{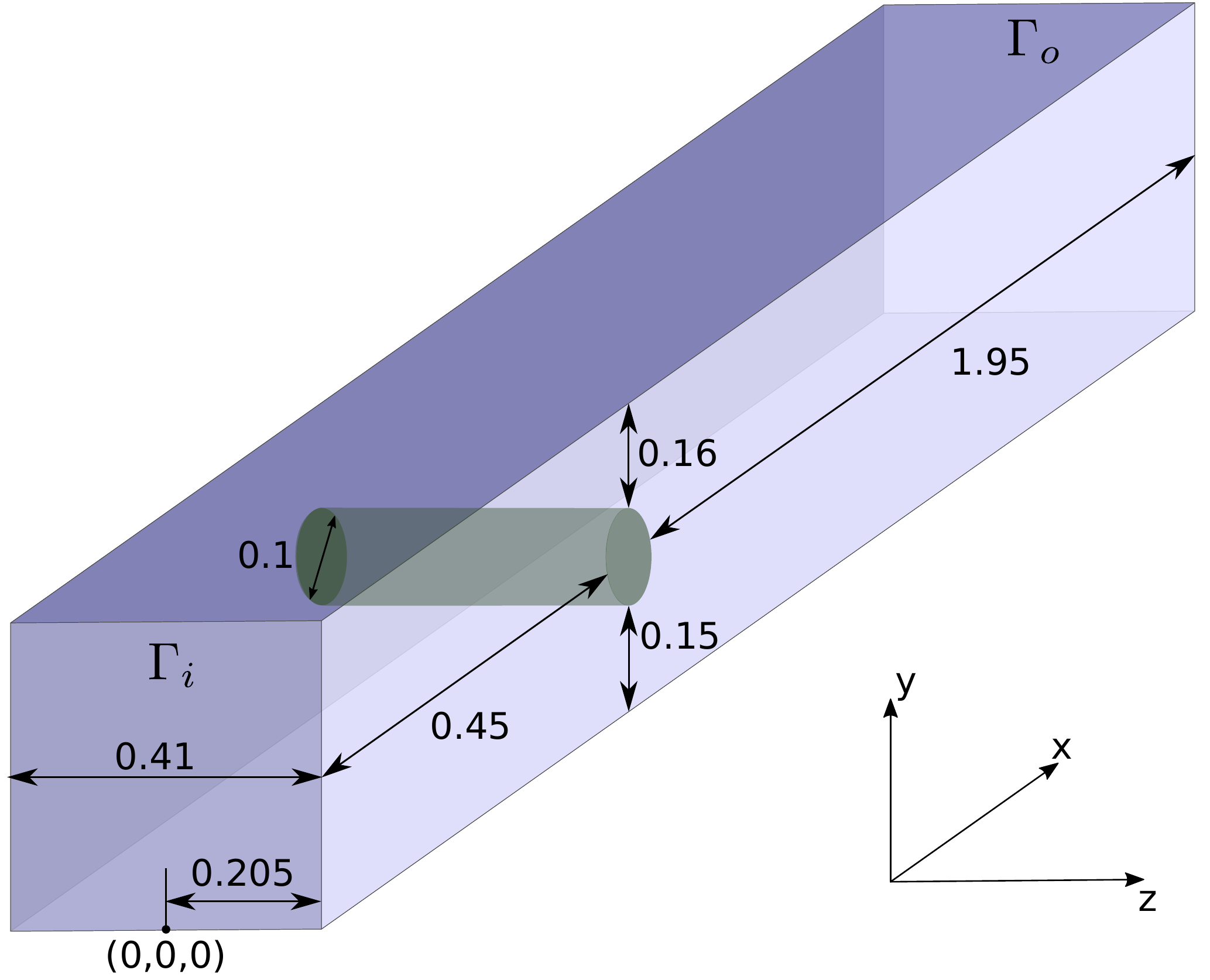}
	\caption{Geometrical setting of the 3d benchmark problem.} 
	\label{Fig:3d_benchmark_setting}
\end{figure}

In this subsection the proposed GMG approach is applied to simulate flow around a cylinder in three space dimension;  cf.\, \cite{johnEfficiencyLinearizationSchemes2006,schaferBenchmarkComputationsLaminar1996}. We note that this benchmark continues to be a challenging test problem for flow solvers.
So far, the benchmark is still an open one since guaranteed numbers for the drag and lift coefficients are not available yet.
The geometry of the benchmark is shown in Fig.~\ref{Fig:3d_benchmark_setting}.

The goal quantities are again the drag and lift coefficients. With the drag and lift forces defined in \eqref{eq:drag_lift_forces}, the drag and lift coefficients are given by
\begin{align}
	c_D &= \frac{2}{\bar{U}^2 D H} F_D \,,
	&
	c_L &= \frac{2}{\bar{U}^2 D H} F_L \,,
\end{align}
where the diameter of the cylinder is $D = 0.1$ and the height of the pipe is $H = 0.41$ (cf.\ Fig.~\ref{Fig:3d_benchmark_setting}). On the inflow boundary $\Gamma_i$ the fluid velocity $\vec{v} = (v_x, v_y, v_y)^{\top}$ with 
\begin{align}
	v_x(\vec{x})
	&=
	\frac{- 16 \cdot U_m ~\cdot y \cdot \left( z - \frac{H}{2} \right)
		\cdot (H - y) \cdot \left(z + \frac{H}{2}\right)}{H^4}\,,
	&
	v_y(\vec{x}) &= 0\,,
	&
	v_z(\vec{x}) &= 0\,,
\end{align}
and $U_m = 2.25$ is prescribed. By the characteristic velocity of the flow of $\overline{U} = 1$ and a viscosity of $\nu = 0.001$ we compute the Reynolds number of the flow to  
\begin{equation}
	Re = \frac{\overline{U} \cdot D}{\nu} = 100 \,.
\end{equation}
The final simulation time is put to $T=8$ such that $I = (0, 8]$.

\begin{figure}[h!tb]
	\centering
	\includegraphics[width=1.0\textwidth,keepaspectratio]
	{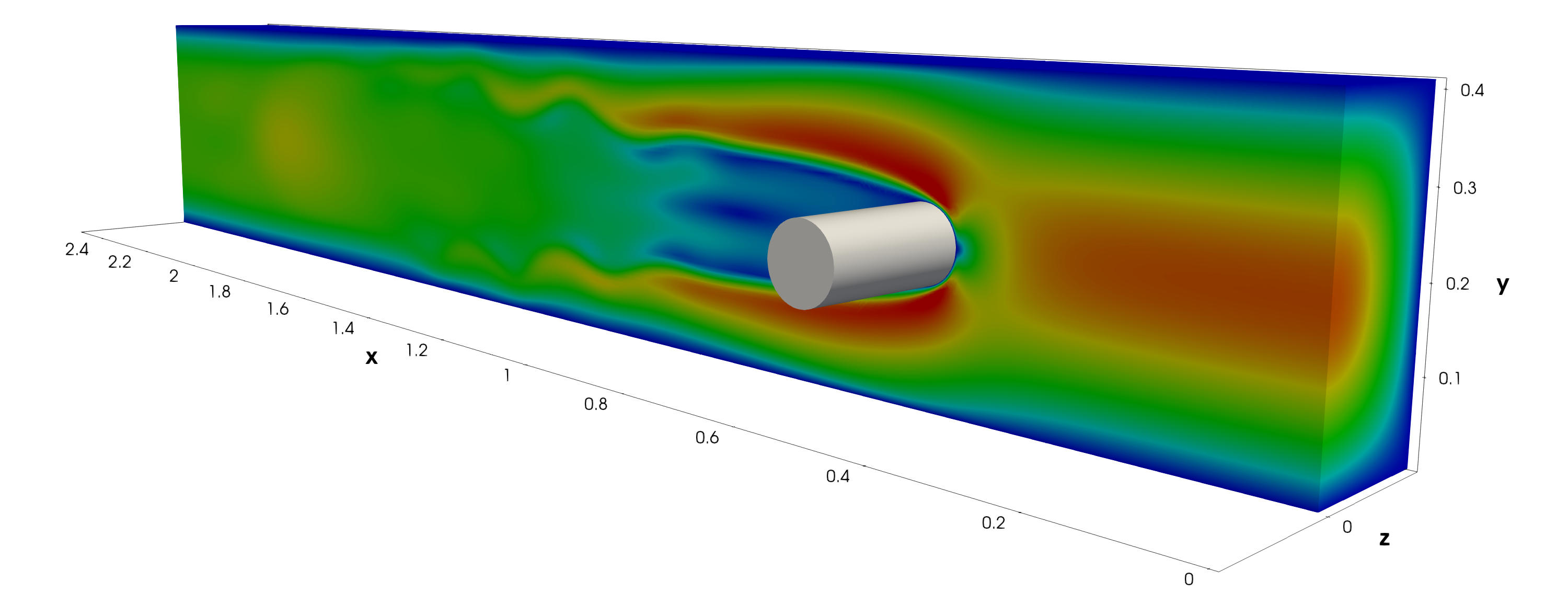}
	\caption{Flow profile in the longitudinally cut domain (at $z = 0$) of the benchmark problem with $Re = 100$.} 
	\label{Fig:3d_benchmark_1}
\end{figure}

The numerical approximation is done in the space-time finite element spaces $(\mathbb P_1(I_n;H_h^2)^2\times \mathbb P_1(I_n;H_{h,\text{disc}}^1)$.
Thus, the discontinuous Galerkin approximation in time with piecewise linear polynomials is used.
\textblue{We perform threes simulations with different spatial mesh sizes, shown in \cref{tab:3d_DFG_results}}.
\textblue{On the largest problem, Nr. 1}, this results in $\num{96876736}$ spatial degrees of freedom \textblue{on the finest mesh level $G$} in each time interval $I_n$ , i.e., over all degrees of freedom in time on $I_n$.
The time interval $I=(0,T]$ is divided into $\num{1598}$ slices of different length, due to the benchmark configuration.
\textblue{The simulation is performed on up to 32 nodes of the Linux cluster (see \cref{tab:3d_DFG_timings})}.
To each CPU core an own process is assigned. Thus, the simulations \textblue{is run e.g. in setting Nr. 1 by $32 \cdot 2 \cdot 14 = 896$ processes.}
\textblue{In setting Nr. 1 the mesh \textblue{level G} consists} of $\num{1703936}$ cells such that each process accesses $\num{1901} \pm 1$ cells of the mesh.
The amount of memory of each process to store all the cell inverses $\mat{J}_K^{-1}$ \textblue{on the finest level in setting Nr. 1} therefore is $\SI{231.2}{\kilo\byte} \cdot 1902 \approx \SI{439.7}{\mega\byte}$. For current high performance computing systems this represents a very reasonable or even small amount of memory usage. \Cref{Fig:3d_benchmark_1} visualizes the computed velocity field in the longitudinally cut domain at the final simulation time $T = 8$.

{
	\sisetup{scientific-notation = false,
		round-mode=places,
		round-precision=4,
		output-exponent-marker=\ensuremath{\mathrm{e}},
		table-figures-integer=1, 
		table-figures-decimal=3, 
		table-figures-exponent=1, 
		table-sign-mantissa = false, 
		table-sign-exponent = true, 
		table-number-alignment=center} 
		
	\begin{table}[h!t]
		\caption{Computed drag and lift coefficients and average number of Newton steps per time step and of GMRES iterations per Newton step in the three-dimensional benchmark \textblue{for three simulations. $h_{\text{max}}$ is the maximum diameter of the cells on the finest mesh level $G$ and DoFs denotes the problem size on this level.}}
		\centering
		\begin{tabular}{c@{\hskip 4ex} c@{\hskip 4ex} b@{\hskip 4ex} c@{\hskip 4ex}  c@{\hskip 4ex} c@{\hskip 4ex} c@{\hskip 4ex} b@{\hskip 4ex}}
	\toprule
	{Nr.} & {DoFs} & {$h_{\text{max}}$} & {$c_{D_{max}}$} & {$c_{L_{max}}$} & {$\bar{n}_{\text{Newton}}$} & {$\bar{n}_{\text{GMRES}}$}  & {$\bar{n}_{\text{GMRES}}$} \\
	\cmidrule(lr){1-3} \cmidrule(r){4-5} \cmidrule(r){6-6} \cmidrule(r){7-7} \cmidrule(r){8-8}
	{1} & \num{96876736} & \num{0.0110} & \num{3.28768} & \num{-0.00722495} & 1.47 & \num{87} & \num{26} \\
	{2} & \num{12293216} & \num{0.0220} & \num{3.20697} & \num{-0.00336773} & 1.42 & \num{74} & \num{22} \\
	{3} & \num{1583152} & \num{0.0440} & \num{2.93094} &  \num{-0.00312805} & 1.44 & \num{72} & \num{24} \\
	\cmidrule(lr){1-3} \cmidrule(r){4-5} \cmidrule(r){6-6} \cmidrule(r){7-7} \cmidrule(r){8-8}
	\multicolumn{6}{c}{Multigrid cycle} & {$V(1,1)$} & {$V(4,4)$} \\
	\bottomrule
\end{tabular}
		\label{tab:3d_DFG_results}
	\end{table}
}

\begin{figure}[h!tb]
	\centering
	\subcaptionbox{Computed drag coefficients $c_D$. \label{fig:drag_coefficients}}
	[1.\columnwidth]
	{\includegraphics[width=0.98\textwidth,keepaspectratio]{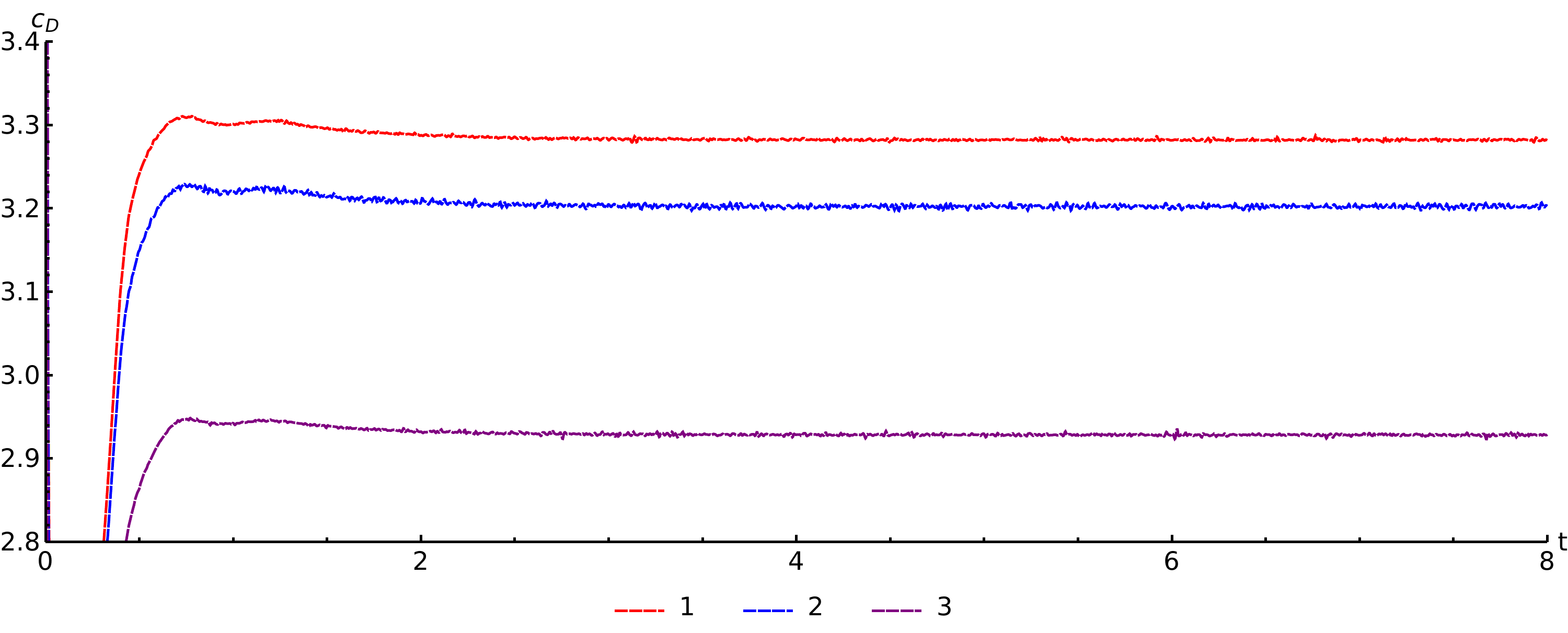}	
	}
	\\[2ex]
	\subcaptionbox{Computed lift coefficients $c_L$. \label{fig:lift_coefficients}}
	[1.\columnwidth]
	{\includegraphics[width=0.98\textwidth,keepaspectratio]{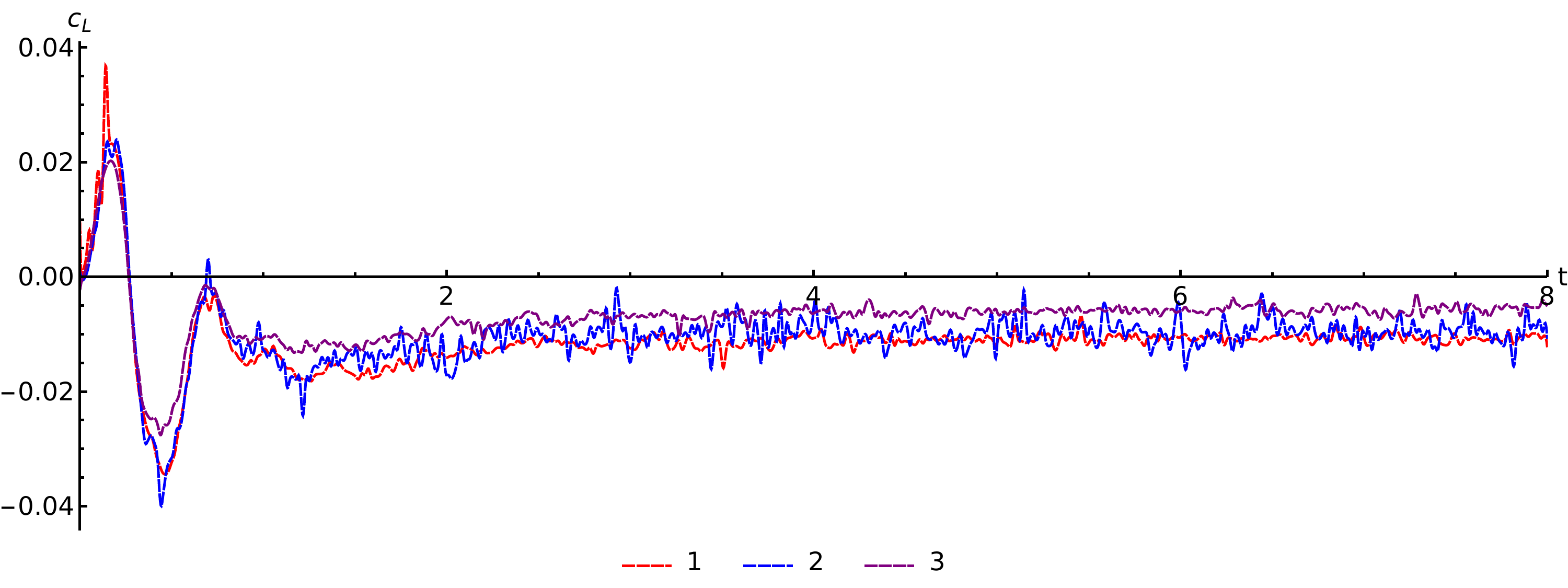}	
	}	
	\caption{Computed drag and lift coefficients of the 3d benchmark on different mesh levels.}
	\label{fig:drag_lift}
\end{figure}

Table~\ref{tab:3d_DFG_results} and Fig.~\ref{fig:drag_lift} present the computed drag and lift coefficients.
\textblue{The coarse level is always set to $g = 1$ and we used 6 - Nr. mesh levels. So for instance for simulation Nr. 1 the finest mesh level is level 5.}
Moreover, Table~\ref{tab:3d_DFG_results} summarizes the average number of Newton steps per subinterval and GMRES iterations per Newton step.

The efficiency of the Newton iteration for solving the nonlinear problem is clearly demonstrated. The average number of Newton iterations is smaller than in the two-dimensional case; cf.\ Table~\ref{tab:2d_DFG_results}. This might be due to fact that the stopping criteria was weakened to a tolerance of $\num{1e-6}$ instead of $\num{1e-8}$ in the two-dimensional case.
Again, the GMG preconditioned GMRES solver shows an almost grid independent convergence behavior. The average number of GMRES iterations per Newton steps is only increased very slightly by grid refinements.
Thereby, the high efficiency of the proposed GMG preconditioning is demonstrated impressively.  

Finally, Table~\ref{tab:3d_DFG_timings} shows the wall-time consumption of the code for three mesh levels of successive refinement in space. In contrast to the two-dimensional case, most of the compute time is now spent on solving the Newton-linearized system. The main reason for this shift is probably the increased number of GMRES steps, compared to Table~\ref{tab:2d_DFG_results}, that are performed until convergence of the GMRES method is reached.

{
	\sisetup{scientific-notation = false,
		round-mode=places,
		round-precision=2,
		output-exponent-marker=\ensuremath{\mathrm{e}},
		table-figures-integer=1, 
		table-figures-decimal=3, 
		table-figures-exponent=1, 
		table-sign-mantissa = false, 
		table-sign-exponent = true, 
		table-number-alignment=center} 
		
	\begin{table}[h!t]
		\small
		\caption{Wall time consumption of the three-dimensional benchmark simulation.}
	\label{tab:3d_DFG_timings}
\begin{subtable}[t]{\textwidth}
	\centering
	\caption{Utilizing a $V(1,1)$ multigrid cycle.}
	\label{tab:3d_DFG_timings_v_1_1}
	\begin{tabular}{c@{\,\,\,\,}cc  c@{\,\,\,\,}c   c@{\,\,\,\,}c  c@{\,\,\,\,}c}
		\toprule
		{DoFs} & {$n_{\text{Nodes}}$} & {$t_{\text{wall}}$} &
		{ $t_{\text{GMRES}}$ } & {\% of $t_{wall}$} &
		{$t_{\text{Inv}}$} & {\% of $t_{wall}$} &
		{$t_{\text{Upd}}$} & {\% of $t_{wall}$} \\
		\cmidrule(lr){1-3}
		\cmidrule(lr){4-5}
		\cmidrule(lr){6-7}
		\cmidrule(lr){8-9}
		\num{96876736} & \num{32} &\SI{153.89}{\hour} &
		\SI{92.22}{\hour} & \num{60.00} &
		\SI{1.57}{\hour} & \num{1.00} &
		\SI{0.51}{\hour} & \num{0.44} \\
		\num{12293216} & \num{16} &\SI{45.50}{\hour} &
		\SI{27.90}{\hour} & \num{65.64} &
		\SI{0.14}{\hour} & \num{0.34} &
		\SI{0.07}{\hour} & \num{0.16} \\
		\num{1583152} & \num{2} & \SI{7.36}{\hour} &
		\SI{5.22}{\hour} & \num{70.94} &
		\SI{0.04}{\hour} & \num{0.74} &
		\SI{0.04}{\hour} & \num{0.50} \\
		\bottomrule
	\end{tabular}
\end{subtable}
\begin{subtable}[t]{\textwidth}
	\vspace*{3ex}
	\centering
	\caption{\textblue{Utilizing a $V(4,4)$ multigrid cycle.}}
	\textblue{
	\label{tab:3d_DFG_timings_v_4_4}
	\begin{tabular}{c@{\,\,\,\,}cc  c@{\,\,\,\,}c   c@{\,\,\,\,}c  c@{\,\,\,\,}c}
		\toprule
		{DoFs} & {$n_{\text{Nodes}}$} & {$t_{\text{wall}}$} &
		{ $t_{\text{GMRES}}$ } & {\% of $t_{wall}$} &
		{$t_{\text{Inv}}$} & {\% of $t_{wall}$} &
		{$t_{\text{Upd}}$} & {\% of $t_{wall}$} \\
		\cmidrule(lr){1-3}
		\cmidrule(lr){4-5}
		\cmidrule(lr){6-7}
		\cmidrule(lr){8-9}
		\num{96876736} & \num{32} &\SI{106.99}{\hour} &
		\SI{45.32}{\hour} & \num{42.36} &
		\SI{1.70}{\hour} & \num{1.59} &
		\SI{0.47}{\hour} & \num{0.44} \\
		\num{12293216} & \num{16} &\SI{32.00}{\hour} &
		\SI{14.40}{\hour} & \num{45.00} &
		\SI{0.15}{\hour} & \num{0.46} &
		\SI{0.08}{\hour} & \num{0.24} \\
		\num{1583152} & \num{2} & \SI{4.96}{\hour} &
		\SI{2.82}{\hour} & \num{56.86} &
		\SI{0.04}{\hour} & \num{0.77} &
		\SI{0.04}{\hour} & \num{0.83} \\
		\bottomrule
	\end{tabular}
	}
\end{subtable}
	\end{table}
}

\section{Summary and outlook}
\label{Sec:SumOut}

In this work a parallel GMG preconditioner with a cell-based Vanka smoother for solving the nonstationary, incompressible Navier--Stokes equations was presented. Its efficient implementation in the deal.II finite element library was discussed. Discontinuous Galerkin methods and inf-sup stable pairs of finite element spaces with discontinuous pressure elements were used for the discretization of the time and space variables, respectively. The GMG preconditioner was applied to a flexible GMRES method for solving the Newton linearized algebraic problem. The performance properties of the GMG method and its parallel implementation were analyzed computationally for the two- and three-dimensional benchmark problem of flow around a cylinder. A quasi grid independence of the GMG preconditioned GMRES solver was observed confirming the high efficiency of the GMG approach. In a forthcoming work we will address an extension of the proposed GMG method to discretizations of the Navier--Stokes equations on evolving domains by using CutFEM techniques on fixed background meshes; cf.~\cite{anselmannCutFiniteElement2022}.

\section*{Acknowledgments}

The authors wish to thank Friedhelm Schieweck from the University of Magdeburg for his helpful support to the development and implementation of the GMG approach.

\bibliographystyle{IEEEtranS}
\bibliography{library} 

\end{document}